\documentclass[12pt]{article}
\usepackage{amssymb}
\usepackage{amsmath}
\usepackage[mathscr]{eucal}

\renewcommand{\phi}{\varphi}

\newcommand{\up}[1]{\stackrel{\raisebox{-0.07cm}{\tiny$\smile$}}{#1}}
\newcommand{\down}[1]{\stackrel{\raisebox{-0.08cm}{\tiny$\frown$}}{#1}}
\begin{document}
\newcounter{abschnitt}
\newtheorem{koro}{Corollary}[abschnitt]
\newtheorem{defi}[koro]{Definition}
\newtheorem{satz}[koro]{Theorem}
\newtheorem{lem}[koro]{Lemma}
\newtheorem{bem}[koro]{Bemerkungen}
\newtheorem{conj}[koro]{Conjecture}

\newcounter{saveeqn}
\newcommand{\alpheqn}{\setcounter{saveeqn}{\value{abschnitt}}%
\renewcommand{\theequation}{%
\mbox{\arabic{saveeqn}.\arabic{equation}}}}%
\newcommand{\reseteqn}{\setcounter{equation}{0}%
\renewcommand{\theequation}{\arabic{equation}}}

\begin{center}
{\Large\bf Convolutions and Multiplier Transformations \\ of
Convex Bodies}
\end{center}

\begin{center}
Franz E. Schuster
\end{center}

\begin{quote}
\footnotesize{ \vskip 1truecm\noindent {\bf Abstract.} Rotation
intertwining maps from the set of convex bodies in $\mathbb{R}^n$
into \linebreak itself that are continuous linear operators with
respect to Minkowski and Blaschke addition are investigated. The
main focus is on Blaschke-Minkowski homomorphisms. We show that
such maps are represented by a spherical convolution operator. An
\linebreak application of this representation is a complete
classification of all even Blaschke-Minkowski homomorphisms which
shows that these maps behave in many respects similar to the well
known projection body operator. Among further applications is the
\linebreak following result: If an even Blaschke-Minkowski
homomorphism maps a convex body to a polytope, then it is a
constant multiple of the projection body operator.

}
\end{quote}

\vspace{1cm}

\centerline{\large{\bf{ \setcounter{abschnitt}{1}
\arabic{abschnitt}. Introduction and Statement of Main Results}}}

\alpheqn

\vspace{0.7cm}\noindent For $n \geq 3$ let $\mathcal{K}^n$ be the
set of convex bodies in $\mathbb{R}^n$, i.e.\ nonempty, compact,
convex sets, equipped with the Hausdorff topology. Let
$\mathcal{K}_0^n$ be the subset of $\mathcal{K}^n$ consisting of
the convex bodies with interior points. A convex body $K \in
\mathcal{K}^n$ is determined by its support function $h(K,\cdot)$,
defined on $\mathbb{R}^n$ by $h(K,x)=\max\{x\cdot y:y \in K\}$.
We shall mostly consider the restriction of support functions to
the Euclidean unit sphere $S^{n-1}$ which are elements of
$\mathcal{C}(S^{n-1})$, the space of continuous functions on
$S^{n-1}$ with the uniform topology.

By Minkowski's existence theorem, a convex body $K \in
\mathcal{K}_0^n$ is also determined up to translation by its
surface area measure (of order $n - 1$) $S_{n-1}(K,\cdot)$. The
measure of a Borel set $\omega \subseteq S^{n-1}$ is the $n-1$
dimensional Hausdorff measure of the set of all boundary points
of $K$ at which there exists a normal vector of $K$ belonging to
$\omega$. $S_{n-1}(K,\cdot)$ is an element of
$\mathcal{M}_o^+(S^{n-1})$, the space of nonnegative measures on
$S^{n-1}$ having their center of mass in the origin, equipped with
the $\mbox{weak}^*$ topology.

Two of the most important algebraic structures on the set of
convex bodies are Minkowski (vector) addition and Blaschke
addition. For $K_1, K_2 \in \mathcal{K}^n$ and $\lambda_1,
\lambda_2 \geq 0$, the support function of the Minkowski linear
combination $\lambda_1K_1+\lambda_2K_2$ is
\[h(\lambda_1K_1+\lambda_2K_2,\cdot)=\lambda_1h(K_1,\cdot)+\lambda_2h(K_2,\cdot).\]

For $K_1, K_2 \in \mathcal{K}_0^n$ and $\lambda_1, \lambda_2 \geq
0$ (not both 0), the Blaschke linear combination $\lambda_1\cdot
K_1 \: \# \: \lambda_2\cdot K_2$ is defined (up to translation) by
\begin{equation*}
S_{n-1}(\lambda_1\cdot K_1 \: \# \: \lambda_2\cdot
K_2,\cdot)=\lambda_1S_{n-1}(K_1,\cdot)+\lambda_2S_{n-1}(K_2,\cdot).
\end{equation*}
With these operations $\mathcal{K}^n$ and $[\mathcal{K}_0^n]$, the
set of translation classes of convex bodies with nonempty
interior, are abelian semi-groups.

\pagebreak

There are several groups acting on the spaces $\mathcal{K}^n$ and
$[\mathcal{K}_0^n]$. One of the most important is the group of
rotations $SO(n)$.

With these different structures on the space of convex bodies it
is natural to ask what kind of operators on $\mathcal{K}^n$ and
$[\mathcal{K}_0^n]$ are compatible with the algebraic and
topological structures. From a geometric point of view we are
interested in operators that intertwine rotations.

In \cite{schneider74} Schneider started an investigation of
continuous rigid motion intertwining and Minkowski additive maps
which he called Minkowski endomorphisms. Among other results he
obtained (under additional assumptions) characterizations of
several such mappings. In \cite{schneider74b} Schneider
classified all Minkowski endomorphisms in $\mathbb{R}^2$.
Kiderlen, relaxing the definition of Minkowski endomorphisms to
continuous {\it rotation} intertwining and Minkowski additive
maps, extended in \cite{kiderlen05} Schneider's classification
result to higher dimensions for weakly monotone Minkowski
endomorphisms, i.e.\ they are monotone with respect to set
inclusion on convex bodies having their Steiner point in the
origin. Kiderlen also gave a complete classification of all
Blaschke endomorphisms, i.e.\ continuous rotation intertwining and
Blaschke additive maps.

\begin{defi} \label{defbmhomo} A map $\Phi: [\mathcal{K}_0^n] \rightarrow
\mathcal{K}^n$ is called Blaschke-Minkowski homomorphism if it
satisfies the following conditions:
\begin{enumerate}
\item $\Phi$ is continuous.
\item For all $K, L \in [\mathcal{K}_0^n]$,
\begin{equation} \label{bmadd}
\Phi(K \: \# \: L)=\Phi K + \Phi L.
\end{equation}
\item $\Phi$ is rotation intertwining, i.e.\ for all $K \in
[\mathcal{K}_0^n]$ and every $\vartheta \in SO(n)$,
\[\Phi(\vartheta K)=\vartheta \Phi K.   \]
\end{enumerate}
\end{defi}

The well known projection body operator, see (\ref{projbody}), is
an example of a Blaschke-Minkowski homomorphism. For its many
applications in different areas see
\cite{bourgainlindenstrauss88}, \cite{bolker69},
\cite{gardner95}, \cite{goodeyweil93}, \cite{grinbergzhang99},
\cite{groemer96} and \cite{ludwig02}. Further examples of
Blaschke-Minkowski homomorphisms can be found in
\cite{goodeyweil92} and \cite{schneiderhug02}. The operator that
maps every convex body to the origin is called the {\it trivial}
Blaschke-Minkowski homomorphism.

The main purpose of this article is to show that there is a
representation for Blaschke-Minkowski homomorphisms analogous to
the ones obtained by Schneider and Kiderlen and to establish a
connection to the theory of Minkowski and Blaschke endomorphisms
developed by them. Moreover we will characterize special Blaschke-
Minkowski homomorphisms and investigate the volume and more
general quermassintegrals of images under these mappings.

Classification results of mappings of convex bodies, in particular
of valuations on convex sets, form a main part of convex
geometry. Here, a map $\Phi$ defined on $\mathcal{K}^n$ and
taking values in an abelian semigroup is called a \emph{valuation}
if for all $K, L \in \mathcal{K}^n$ such that also $K \cup L \in
\mathcal{K}^n$,
\[\Phi(K \cup L)+\Phi(K \cap L)=\Phi K + \Phi L.   \]
The theory of valuations and its important applications in
integral geometry and geometric probability are developed and
described in \cite{hadwiger51}, \cite{klainrota97},
\cite{mcmullen93}, \cite{mcmullenschneider83}.

In \cite{ludwig02} and \cite{ludwig04} valuations taking values
in $\mathcal{K}^n$ are investigated which intertwine volume
preserving linear transformations. In Section 4 we will see that
Blaschke-Minkowski homomorphisms are special valuations. Thus,
our results are a contribution to the classification of
continuous rotation intertwining valuations taking values in
$\mathcal{K}^n$. In full generality this problem is still open.

Let $SO(n)$ denote the group of rotations in $n$ dimensions.
Identifying $S^{n-1}$ with the homogeneous space $SO(n)/SO(n-1)$,
where $SO(n-1)$ denotes the group of rotations leaving the point
$\mbox{\raisebox{-0.01cm}{$\down{e}$}}$ (the pole) of $S^{n-1}$
fixed, it is possible to introduce a natural convolution
structure on $\mathcal{C}(S^{n-1})$ and $\mathcal{M}(S^{n-1})$. A
special role play convolution operators generated by $SO(n-1)$
invariant (or zonal) functions and measures. The set of continuous
zonal functions on $S^{n-1}$ will be denoted by
$\mathcal{C}(S^{n-1},\mbox{\raisebox{-0.01cm}{$\down{e}$}})$.

A function $f \in \mathcal{C}(S^{n-1})$ is called weakly positive
if there exists a vector $x \in \mathbb{R}^n$ such that $f(u)+x
\cdot u \geq 0$ for every $u \in S^{n-1}$. The main theorem of
this article is the following representation for
Blaschke-Minkowski homomorphisms:

\begin{satz} \label{satzbmhomo} If $\Phi: \mathcal{K}^n \rightarrow \mathcal{K}^n$ is a
Blaschke-Minkowski homomorphism, then there is a weakly positive
$g \in
\mathcal{C}(S^{n-1},\mbox{\raisebox{-0.01cm}{$\down{e}$}})$,
unique up to addition of a linear function, such that
\begin{equation} \label{bmhomorep}
h(\Phi K,\cdot)=S_{n-1}(K,\cdot) \ast g.
\end{equation}
\end{satz}

Note that in Theorem \ref{satzbmhomo} the domain of $\Phi$ is the
set $\mathcal{K}^n$ in contrast to Definition \ref{defbmhomo}. The
reason for this is a natural identification of maps on
$[\mathcal{K}_0^n]$ with translation invariant maps on
$\mathcal{K}_0^n$ and the fact (as we will show) that there is a
unique continuous extension of every Blaschke-Minkowski
homomorphism to $\mathcal{K}^n$.

A map $\Phi: \mathcal{K}^n \rightarrow \mathcal{K}^n$ is called
\emph{even} if $\Phi K=\Phi(-K)$ for every $K \in \mathcal{K}^n$.
We call a body $K \in \mathcal{K}^n$ a body of revolution if $K$
is invariant under rotations of $SO(n-1)$. Using Theorem
\ref{satzbmhomo} and a further investigation of properties of
generating functions of Blaschke-Minkowski homomorphisms, a
classification of all even Blaschke-Minkowski homomorphisms is
possible.

\begin{satz} \label{bmsymm} A map $\Phi: \mathcal{K}^n \rightarrow \mathcal{K}^n$
is an even Blaschke-Minkowski homomorphism \linebreak if and only
if there is a centrally symmetric body of revolution $L \in
\mathcal{K}^n$, unique up to translation, such that
\[h(\Phi K,\cdot)=S_{n-1}(K,\cdot) \ast h(L,\cdot).\]
\end{satz}

The projection body operator $\Pi: \mathcal{K}^n \rightarrow
\mathcal{K}^n$ is defined by
\begin{equation} \label{projbody}
h(\Pi
K,u)=\mbox{vol}_{n-1}(K|u^{\bot})=\frac{1}{2}(S_{n-1}(K,\cdot)\ast
h([-\mbox{\raisebox{-0.01cm}{$\down{e}$}},\mbox{\raisebox{-0.01cm}{$\down{e}$}}],\cdot))(u),
\end{equation}
where
$[-\mbox{\raisebox{-0.01cm}{$\down{e}$}},\mbox{\raisebox{-0.01cm}{$\down{e}$}}]$
denotes the segment with endpoints
$-\mbox{\raisebox{-0.01cm}{$\down{e}$}}$ and
$\mbox{\raisebox{-0.01cm}{$\down{e}$}}$. The operator $\Pi$ maps
polytopes to finite Minkowski linear combinations of rotated and
dilated copies of the line segment
$[-\mbox{\raisebox{-0.01cm}{$\down{e}$}},\mbox{\raisebox{-0.01cm}{$\down{e}$}}]$,
which is a geometric interpretation of the convolution formula in
(\ref{projbody}). A general convex body is mapped by $\Pi$ to a
zonoid, i.e.\ a limit of Minkowski sums of line segments. By
Theorem \ref{bmsymm}, a general even Blaschke-Minkowski
homomorphism maps polytopes to finite Minkowski linear
combinations of rotated and dilated copies of a symmetric body of
revolution $L$. General convex bodies are again mapped to limits
of these finite Minkowski linear combinations.

\pagebreak

In \cite{kiderlen05} a notion of adjointness between Minkowski and
Blaschke endomorphisms was introduced. The following consequence
of Theorem \ref{satzbmhomo} illustrates the behaviour of adjoint
endomorphisms in conjunction with Blaschke-Minkowski
homomorphisms.

\begin{satz} \label{bmanddual} Let $\Psi$ be a Minkowski and $\Psi^*$ a Blaschke endomorphism.
Then the following statements are equivalent:
\begin{enumerate}
\item $\Psi$ and $\Psi^*$ are adjoint endomorphisms.
\item For every Blaschke-Minkowski homomorphism $\Phi$
\begin{equation} \label{phipsipsiphi}
\Phi \circ \Psi^* = \Psi \circ \Phi.
\end{equation}
\item (\ref{phipsipsiphi}) holds for some injective Blaschke-Minkowski homomorphism
$\Phi$.
\end{enumerate}
\end{satz}

A different application of Theorem \ref{bmsymm} is the following
characterization of $\Pi$.

\begin{satz} \label{bmsuppnopoly} Let $\Phi: \mathcal{K}^n \rightarrow \mathcal{K}^n$
be an even Blaschke-Minkowski homomorphism. If there exists a
convex body $K \in \mathcal{K}_0^n$ such that $\Phi K$ is a
polytope, then there is a constant $c \in \mathbb{R}^+$ such that
\[\Phi = c \Pi.   \]
\end{satz}

As a consequence of Theorem \ref{satzbmhomo} the image of a
Minkowski linear combination under a Blaschke-Minkowski
homomorphism is a homogeneous polynomial of degree $n-1$. In
particular, Blaschke-Minkowski homomorphisms satisfy, for $K \in
\mathcal{K}^n$, the Steiner type formula
\[\Phi(K+\varepsilon B^n)=\sum \limits_{i=0}^{n-1}\varepsilon^{i} {n-1 \choose i}\Phi_{i}K,   \]
where $B^n$ is the Euclidean unit ball and the sum is with
respect to Minkowski addition. The operators $\Phi_i:
\mathcal{K}^n \rightarrow \mathcal{K}^n$, $i=0,\ldots,n-1$, are
continuous and rotation intertwining. The image of a ball under a
Blaschke-Minkowski homomorphism $\Phi$ is  again a ball. Let in
the following $r_\Phi \in \mathbb{R}^+$ denote the radius of
$\Phi B^n$ and $\kappa_n$ the volume of $B^n$.

We will prove a strengthened version of the classical inequality
between the two consecutive quermassintegrals $W_{n-1}$ and
$W_{n-2}$, using the induced weakly monotone Minkowski
endomorphisms $\Phi_{n-2}$.

\begin{satz} \label{minkinequ} Let $\Phi: \mathcal{K}^n \rightarrow \mathcal{K}^n$ be a nontrivial Blaschke-Minkowski homomorphism.
If $K \in \mathcal{K}^n$, then
\begin{equation} \label{sharpmin}
W_{n-1}(K)^2 \geq \frac{\kappa_n}{r_{\Phi}^2}W_{n-2}(\Phi_{n-2}K)
\geq \kappa_nW_{n-2}(K).
\end{equation}
If $K$ is not a singleton, there is equality on the left hand side
only if $\Phi_{n-2}K$ is a ball and equality on the right hand
side only if $K$ is ball.
\end{satz}

Inequality (\ref{sharpmin}) is related to a conjectured
projection inequality of Petty for the volume of projection
bodies, see \cite{gardner95}, \cite{lutwak90b} and \cite{petty72}.

\pagebreak

\setcounter{abschnitt}{2}
\centerline{\large{\bf{\arabic{abschnitt}. Spherical Convolution
and Spherical Harmonics}}}

\reseteqn \alpheqn

\setcounter{koro}{0}

\vspace{0.7cm} \noindent As we deal with different kinds of
analytical representations of convex bodies by functions and
measures on $S^{n-1}\cong SO(n)/SO(n-1)$, we will first introduce
some basic notions connected to $SO(n)$ and $S^{n-1}$. As general
reference for this section we recommend the article by Grinberg
and Zhang \cite{grinbergzhang99} and the book by Groemer
\cite{groemer96}.

The identification of $S^{n-1}$ with $SO(n)/SO(n-1)$ is for $u
\in S^{n-1}$ given by
\[u=\vartheta\mbox{\raisebox{-0.01cm}{$\down{e}$}} \mapsto \vartheta SO(n-1).\]
The projection from $SO(n)$ onto $S^{n-1}$ is $\vartheta \mapsto
\mbox{\raisebox{0.04cm}{$\down{\vartheta}$}}:=\vartheta
\mbox{\raisebox{-0.01cm}{$\down{e}$}}$. The unity $e \in SO(n)$ is
mapped to the pole of the sphere
$\mbox{\raisebox{-0.01cm}{$\down{e}$}} \in S^{n-1}$. $SO(n)$ and
$S^{n-1}$ will be equipped with the invariant probability
measures denoted by $d\vartheta$ and $du$.

Let $\mathcal{C}(SO(n))$ denote the set of continuous functions
on $SO(n)$ with the uniform topology and $\mathcal{M}(SO(n))$ its
dual space of signed finite measures on $SO(n)$ with the
$\mbox{weak}^*$ topology. Let $\mathcal{M}^+(SO(n))$ be the set
of nonnegative measures on $SO(n)$. For $\mu \in
\mathcal{M}(SO(n))$ and $f \in \mathcal{C}(SO(n))$, the canonical
pairing is
\begin{equation*}
\langle \mu,f\rangle=\langle f,\mu
\rangle=\int_{SO(n)}f(\vartheta)d\mu(\vartheta).
\end{equation*}
Sometimes we will identify a continuous function $f$ with the
absolute continuous measure with density $f$ and thus view
$\mathcal{C}(SO(n))$ as a subspace of $\mathcal{M}(SO(n))$. The
canonical pairing is then consistent with the usual inner product
on $\mathcal{C}(SO(n))$.

For $\vartheta \in SO(n)$, the left translation $\vartheta f$ of
$f \in \mathcal{C}(SO(n))$ is defined by
\begin{equation} \label{rotate1}
\vartheta f(\eta)=f(\vartheta^{-1}\eta).
\end{equation}
For $\mu \in \mathcal{M}(SO(n))$, we set
\begin{equation} \label{rotate2}
\langle \vartheta \mu,f\rangle=\langle \mu,\vartheta^{-1}f\rangle,
\end{equation}
then $\vartheta \mu$ is just the image measure of $\mu$ under the
rotation $\vartheta$. For $f \in \mathcal{C}(SO(n))$, the function
$\hat{f} \in \mathcal{C}(SO(n))$ is defined by
\begin{equation} \label{hat1}
\hat{f}(\vartheta)=f(\vartheta^{-1}).
\end{equation}
For a measure $\mu \in \mathcal{M}(SO(n))$, we set
\begin{equation} \label{hat2}
\langle \hat{\mu},f\rangle=\langle \mu,\hat{f}\rangle.
\end{equation}

As $SO(n)$ is a compact Lie group the space $\mathcal{C}(SO(n))$
carries a natural convolution structure. For $f, g \in
\mathcal{C}(SO(n))$, the convolution $f \ast g \in
\mathcal{C}(SO(n))$ is defined by
\[(f \ast g)(\eta)=\int_{SO(n)} f(\eta \vartheta^{-1})g(\vartheta)d\vartheta
=\int_{SO(n)}f(\vartheta)g(\vartheta^{-1}\eta)d\vartheta.\]

For $\mu \in \mathcal{M}(SO(n))$, the convolutions $\mu \ast f \in
\mathcal{C}(SO(n))$ and $f \ast \mu \in \mathcal{C}(SO(n))$ with
a function $f \in \mathcal{C}(SO(n))$ are defined by
\begin{equation}\label{convwithmeasure}
(f \ast \mu)(\eta)=\int_{SO(n)} f(\eta
\vartheta^{-1})d\mu(\vartheta), \quad (\mu \ast
f)(\eta)=\int_{SO(n)}\vartheta f(\eta)d\mu(\vartheta).
\end{equation}
Using (\ref{convwithmeasure}), one easily checks that for $\sigma
\in \mathcal{M}(SO(n))$ and $f,g \in \mathcal{C}(SO(n))$
\begin{equation} \label{adjoint}
\langle g \ast \sigma ,f \rangle = \langle g,f \ast
\hat{\sigma}\rangle.
\end{equation}
This leads to the definition of the convolution of two measures
$\mu, \sigma \in \mathcal{M}(SO(n))$
\begin{equation} \label{durchmithut}
\langle\mu \ast \sigma,f \rangle=\langle \sigma,\hat{\mu} \ast
f\rangle = \langle \mu,f\ast \hat{\sigma}\rangle.
\end{equation}

\noindent The convolution on $\mathcal{M}(SO(n))$ is associative,
but as for $n \geq 3$ the group of rotations is not abelian, the
convolution on $\mathcal{M}(SO(n))$ is not commutative. For the
following Lemma see \cite{grinbergzhang99}, p.85.

\begin{lem} \label{approxlem} Let $\mu_m, \mu \in \mathcal{M}(SO(n))$,
$m=1,2,\ldots$ and let $f \in \mathcal{C}(SO(n))$. If $\mu_m
\rightarrow \mu$ weakly, then $f \ast \mu_m \rightarrow f \ast
\mu$ and $\mu_m \ast f \rightarrow \mu \ast f$ uniformly.
\end{lem}

In order to define a convolution structure on
$\mathcal{C}(S^{n-1})$, we will use the method from Grinberg and
Zhang \cite{grinbergzhang99} identifying $S^{n-1}$ with
$SO(n)/SO(n-1)$. This leads to the identification of
$\mathcal{C}(S^{n-1})$ with right $SO(n-1)$-invariant functions
in $\mathcal{C}(SO(n))$ by
\begin{equation} \label{rightinvariant}
\mbox{$\up{f}$}(\vartheta)=f(\vartheta
\mbox{\raisebox{-0.01cm}{$\down{e}$}}),\qquad f \in
\mathcal{C}(S^{n-1}).
\end{equation}
Conversely, every $f \in \mathcal{C}(SO(n))$ induces a continuous
function $\down{f}$ on $S^{n-1}$, defined by
\[\mbox{$\down{f}$}(\mbox{\raisebox{-0.06cm}{$\down{\eta}$}})=\int_{SO(n-1)}f(\eta\vartheta)d\vartheta. \]

\vspace{-0.2cm}

\noindent If $f \in \mathcal{C}(SO(n))$ is right $SO(n-1)$
invariant and $g \in \mathcal{C}(S^{n-1})$ then
$f=\mbox{$\up{\down{f}}$}$ and
$g=\mbox{\raisebox{-0.05cm}{$\down{\up{g}}$}}.$ Thus
$\mathcal{C}(S^{n-1})$ is isomorphic to the subspace of right
$SO(n-1)$ invariant functions in $\mathcal{C}(SO(n))$. For a
measure $\mu \in \mathcal{M}(S^{n-1})$ and a function $f \in
\mathcal{C}(SO(n))$, we set
\[\langle \mbox{\raisebox{-0.05cm}{$\up{\mu}$}},f\rangle =\langle \mu,\down{f}\rangle.\]
In this way the one-to-one correspondence of functions on
$S^{n-1}$ with right $SO(n-1)$ invariant functions on $SO(n)$
carries over to the space $\mathcal{M}(S^{n-1})$ and right
$SO(n-1)$ invariant measures in $\mathcal{M}(SO(n))$.

Note that definitions (\ref{rotate1}), (\ref{rotate2}) and
(\ref{hat1}), (\ref{hat2}) become now meaningful for spherical
functions and measures. Convolution on $\mathcal{C}(S^{n-1}) $ can
be defined via the identification (\ref{rightinvariant}). For
example the convolution of a function $f \in
\mathcal{C}(S^{n-1})$ with a measure $\mu \in
\mathcal{M}(S^{n-1})$ is given by
\[(f \ast \mu)(\mbox{\raisebox{-0.06cm}{$\down{\eta}$}})=(\up{f} \ast \;\mbox{\raisebox{-0.06cm}{$\up{\mu}$}})(\eta)
=\int_{SO(n)}f(\eta\vartheta^{-1}\mbox{\raisebox{-0.02cm}{$\down{e}$}})d\mbox{\raisebox{-0.06cm}{$\up{\mu}$}}(\vartheta).\]
In an analogous way, convolutions of functions or measures can be
defined. Note that \linebreak the Dirac measure
$\delta_{\down{e}}$ is the unique rightneutral element for the
convolution on $S^{n-1}$.

An essential role among spherical functions play $SO(n-1)$
invariant functions. Such a function with the property that
$\vartheta f = f$ for every $\vartheta \in SO(n-1)$, is called
zonal. Zonal functions depend only on the distance of $u$ to
$\mbox{\raisebox{-0.01cm}{$\down{e}$}}$, i.e.\ on the value
$u\cdot \mbox{\raisebox{-0.01cm}{$\down{e}$}}$.

Of course the notion of $SO(n-1)$ invariance carries over to
measures as well. We call a measure $\mu \in \mathcal{M}(S^{n-1})$
zonal, if $\vartheta \mu =\mu$ for every $\vartheta \in SO(n-1)$.
The set of all continuous, zonal functions will be denoted by
$\mathcal{C}(S^{n-1},\mbox{\raisebox{-0.01cm}{$\down{e}$}})$ and
$\mathcal{M}(S^{n-1},\mbox{\raisebox{-0.01cm}{$\down{e}$}})$
denotes the set of zonal measures on $S^{n-1}$.

Spherical convolution becomes simpler for zonal measures. For $f
\in \mathcal{C}(S^{n-1})$ and $\mu \in
\mathcal{M}(S^{n-1},\mbox{\raisebox{-0.01cm}{$\down{e}$}})$, we
have
\begin{equation} \label{zonalconv}
(f\ast \mu)(\mbox{\raisebox{-0.06cm}{$\down{\eta}$}})= \langle
f,\eta \mu \rangle
 =\int_{S^{n-1}}f(\eta u)d\mu(u).
\end{equation}

\noindent For $f \in \mathcal{C}(S^{n-1})$, the rotational
symmetrization $\bar{f} \in
\mathcal{C}(S^{n-1},\mbox{\raisebox{-0.01cm}{$\down{e}$}})$ is
defined by
\[\bar{f}=\delta_{\down{e}}\ast f  =\int_{SO(n-1)}\vartheta fd\vartheta.\]
Since $\delta_{\down{e}}$ is the right invariant element for the
convolution on $S^{n-1}$, we get
\begin{equation} \label{convright}
f \ast g = f \ast \delta_{\down{e}} \ast g = f \ast \bar{g}.
\end{equation}
Thus, for spherical convolution from the right, it suffices to
consider zonal functions and measures. Note that, if $\mu \in
\mathcal{M}(S^{n-1},\mbox{\raisebox{-0.01cm}{$\down{e}$}})$, then
by (\ref{zonalconv}) for every $f \in \mathcal{C}(S^{n-1})$
\begin{equation} \label{rightconvintrot}
(\vartheta f) \ast \mu = \vartheta(f \ast \mu)
\end{equation}
for every $\vartheta \in SO(n)$. Thus the spherical convolution
from the right is a rotation intertwining operator on
$\mathcal{C}(S^{n-1})$ and $\mathcal{M}(S^{n-1})$.

As a zonal function on $S^{n-1}$ depends only on the value of $u
\cdot \mbox{\raisebox{-0.01cm}{$\down{e}$}}$, there is a natural
isomorphism between functions and measures on $[-1,1]$ and zonal
functions and measures on $S^{n-1}$. Define a map $\Lambda:
\mathcal{C}(S^{n-1},\mbox{\raisebox{-0.01cm}{$\down{e}$}})
\rightarrow \mathcal{C}([-1,1]), f \mapsto \Lambda f,$ by
\begin{equation} \label{biglambda}
\Lambda f(t) = f(t
\mbox{\raisebox{-0.01cm}{$\down{e}$}}+\sqrt{1-t^2}v), \qquad v
\in \mbox{\raisebox{-0.01cm}{$\down{e}$}}^\bot \cap S^{n-1}.
\end{equation}
 Then it is easy to see that $\Lambda$ is an
isomorphism with inverse \[\Lambda^{-1}: \mathcal{C}([-1,1])
\rightarrow
\mathcal{C}(S^{n-1},\mbox{\raisebox{-0.01cm}{$\down{e}$}}), f
\mapsto f(\mbox{\raisebox{-0.01cm}{$\down{e}$}}\cdot .\:).\] For
a zonal measure $\mu \in
\mathcal{M}(S^{n-1},\mbox{\raisebox{-0.01cm}{$\down{e}$}})$ and a
function $f \in \mathcal{C}([-1,1])$, define
\[\langle\Lambda \mu,f\rangle=\langle \mu,\Lambda^{-1}f\rangle. \]
The map $\Lambda:
\mathcal{M}(S^{n-1},\mbox{\raisebox{-0.01cm}{$\down{e}$}})\rightarrow
\mathcal{M}([-1,1])$ is the extension of the map defined in
(\ref{biglambda}) and it is again an isomorphism between
$\mathcal{M}(S^{n-1},\mbox{\raisebox{-0.01cm}{$\down{e}$}})$ and
$\mathcal{M}([-1,1])$ with inverse
\[\langle \Lambda^{-1}\mu,f\rangle=\langle \mu,\Lambda \bar{f}\rangle,
\qquad \mu \in \mathcal{M}([-1,1]), f \in \mathcal{C}(S^{n-1}).\]
The isomorphism $\Lambda$ allows one to identify the dual space of
$\mathcal{C}(S^{n-1},\mbox{\raisebox{-0.01cm}{$\down{e}$}})$ with
the space
$\mathcal{M}(S^{n-1},\mbox{\raisebox{-0.01cm}{$\down{e}$}})$.
Using this identification, we obtain for $\mu,\nu \in
\mathcal{M}(S^{n-1},\mbox{\raisebox{-0.01cm}{$\down{e}$}})$ and
$f \in
\mathcal{C}(S^{n-1},\mbox{\raisebox{-0.01cm}{$\down{e}$}})$,
\begin{equation} \label{concom}
\langle \mu \ast \nu,f\rangle=\int_{S^{n-1}}\int_{S^{n-1}}\Lambda
f(u\cdot v)d\mu(u)d\nu(v) =\langle \nu \ast \mu,f\rangle.
\end{equation}
Thus, the convolution of zonal functions and measures is abelian
and $\mathcal{M}(S^{n-1},\mbox{\raisebox{-0.01cm}{$\down{e}$}})$
with the convolution structure becomes an abelian Banach algebra.

\pagebreak

Another property of zonal measures $\mu \in
\mathcal{M}(S^{n-1},\mbox{\raisebox{-0.01cm}{$\down{e}$}})$ is
\begin{equation} \label{selfadjoint}
\hat{\mu}=\mu.
\end{equation}
As a consequence of (\ref{adjoint}) and (\ref{selfadjoint}) we
obtain the following important Lemma.

\begin{lem} \label{selfadlemma}
Let $\mu, \nu \in \mathcal{M}(S^{n-1})$ and $f\in
\mathcal{C}(S^{n-1})$, then
\[\langle \mu \ast \nu ,f \rangle =\langle \mu,f \ast \nu\rangle.\]
\end{lem}

\noindent Using Lemma \ref{selfadlemma} and (\ref{concom}), we get
for $\mu \in \mathcal{M}(S^{n-1})$ and $f \in
\mathcal{C}(S^{n-1},\mbox{\raisebox{-0.01cm}{$\down{e}$}})$,
\begin{equation} \label{zonalconv2}
(\mu \ast f)(u)=\int_{S^{n-1}} \Lambda f(u\cdot v)d\mu(v).
\end{equation}

We will frequently use zonal approximate identities $(\phi_k)_{k
\in \mathbb{N}} $. These are non-negative functions in
$\mathcal{C}^{\infty}(S^{n-1})$. They have already been considered
by Berg \cite{berg69} and we just briefly recall their most
important properties in the following Lemma.

\begin{lem} \label{zonalapproxid} Let $(\phi_k)_{k \in \mathbb{N}}$ be a zonal approximate identity. Then
\begin{enumerate}
\item $f \ast \phi_k \in \mathcal{C}^{\infty}(S^{n-1})$ and $\lim_{k \rightarrow \infty} f \ast \phi_k = f$ uniformly
for every $f \in \mathcal{C}(S^{n-1})$.
\item $\mu \ast \phi_k \in \mathcal{C}^{\infty}(S^{n-1})$ and $\lim_{k \rightarrow \infty} \mu \ast \phi_k = \mu$
weakly for every $\mu \in \mathcal{M}(S^{n-1})$.
\end{enumerate}
\end{lem}

We now collect some facts from the theory of spherical harmonics.
A spherical harmonic of dimension $n$ and order $k$ is the
restriction to $S^{n-1}$ of a harmonic polynomial of order $k$ in
$n$ variables. Let $\mathcal{H}^n_k$ denote the space of
spherical harmonics of dimension $n$ and order $k$.
$\mathcal{H}^n$ will denote the space of all finite sums of
spherical harmonics of dimension $n$.

$\mathcal{H}^n_k$ is a finite dimensional vector space of
dimension $N(n,k)$. The spaces $\mathcal{H}^n_k$ are pairwise
orthogonal with respect to the usual inner product on
$\mathcal{C}(S^{n-1})$. By definition, $\mathcal{H}^n_k$ is
invariant with respect to rotations. Moreover, $\mathcal{H}_k^n$
is irreducible, i.e.\ $\{0\}$ and $\mathcal{H}^n_k$ are the only
subspaces invariant under $SO(n)$. As a consequence we have the
following version of Schur's Lemma for spherical harmonics.

\begin{lem} \label{schur} Let $\Phi: \mathcal{H}_k^n \rightarrow
\mathcal{M}(S^{n-1})$ be a linear map that intertwines rotations.
Then $\Phi$ is either injective or the zero map.
\end{lem}

If $H_1, \ldots, H_{N(n,k)}$ is an orthonormal basis of
$\mathcal{H}_k^n$, then there is a unique polynomial $P_k^n \in
\mathcal{C}([-1,1])$ of degree $k$ such that
\begin{equation} \label{additiontheorem}
\sum \limits_{i=1}^{N(n,k)} H_i(u)H_i(v)=N(n,k)P_k^n(u\cdot v).
\end{equation}
The polynomial $P_k^n$ is called the Legendre polynomial of
dimension $n$ and order $k$. The zonal function $u \mapsto
P_k^n(\mbox{\raisebox{-0.01cm}{$\down{e}$}}\cdot u)$ is up to a
multiplicative constant the unique zonal spherical harmonic in
$\mathcal{H}_k^n$.

The collection $\{H_1, \ldots, H_{N(n,k)}: k \in \mathbb{N}\}$
forms a complete orthogonal system in $\mathcal{L}^2(S^{n-1})$,
i.e.\ for every square integrable function $f$ the series
\[f \sim \sum \limits_{k=0}^{\infty}\pi_k f  \]
converges in quadratic mean to $f$, where $\pi_k f \in
\mathcal{H}_k^n$ is the orthogonal projection of $f$ on the space
$\mathcal{H}_k^n$. Using (\ref{additiontheorem}) and
(\ref{zonalconv}), we obtain
\begin{equation} \label{projectiononhkn}
\pi_k f=\sum\limits_{i=1}^{N(n,k)}\langle f,H_i \rangle H_i =
N(n,k) (f \ast P_k^n(\mbox{\raisebox{-0.01cm}{$\down{e}$}}\cdot
.\:)).
\end{equation}
This leads to the definition of the spherical expansion of a
measure $\mu \in \mathcal{M}(S^{n-1})$
\begin{equation} \label{seriesexmu}
\mu \sim \sum \limits_{k=0}^{\infty}\pi_k\mu,
\end{equation}
where $\pi_k\mu \in \mathcal{H}_k^n$ is defined by
\begin{equation} \label{defpik}
\pi_k\mu=N(n,k) (\mu \ast
P_k^n(\mbox{\raisebox{-0.01cm}{$\down{e}$}}\cdot .\:)).
\end{equation}
We note here two special cases of (\ref{defpik})
\begin{equation} \label{pi0pi1}
\pi_0 \mu = \mu \ast 1 \qquad \mbox{and} \qquad \pi_1 \mu = n \mu
\ast (\mbox{\raisebox{-0.01cm}{$\down{e}$}}\cdot .\:).
\end{equation}
By Lemma \ref{selfadlemma}, we have for every $f \in
\mathcal{C}(S^{n-1})$
\[\langle \pi_k \mu,f\rangle = N(n,k)\langle \mu \ast P_k^n(\mbox{\raisebox{-0.01cm}{$\down{e}$}}\cdot
.\:),f \rangle = N(n,k)\langle \mu,f \ast
P_k^n(\mbox{\raisebox{-0.01cm}{$\down{e}$}}\cdot .\:)\rangle
=\langle \mu,\pi_k f \rangle,\] which, by the completeness of the
system of spherical harmonics, immediately gives:

\begin{lem} \label{uniquedet} Let $\mu \in \mathcal{M}(S^{n-1})$. If
$\mu \ast P_k^n(\mbox{\raisebox{-0.01cm}{$\down{e}$}}\cdot .\:)=0$
for every $k \in \mathbb{N}$ then $\mu = 0$.
\end{lem}

By Lemma \ref{uniquedet}, $\mu \in \mathcal{M}(S^{n-1})$ is
uniquely determined by its series expansion (\ref{seriesexmu}).
Zonal functions and measures are even determined by a sequence of
real numbers. To see this, note that
\[\delta_{\down{e}}\ast P_k^n(u\cdot .\:)=P_k^n(\mbox{\raisebox{-0.01cm}{$\down{e}$}}\cdot u)P_k^n(\mbox{\raisebox{-0.01cm}{$\down{e}$}}\cdot .\:) \]
and thus, by (\ref{zonalconv2}) and (\ref{durchmithut}),
\[(\mu \ast P_k^n(\mbox{\raisebox{-0.01cm}{$\down{e}$}}\cdot .\:))(u)=\langle \mu,P_k^n(u\cdot .\:)\rangle
=\langle \mu,\delta_{\down{e}}\ast P_k^n(u\cdot .\:) \rangle
=\langle \mu,P_k^n(\mbox{\raisebox{-0.01cm}{$\down{e}$}}\cdot
.\:)\rangle P_k^n(\mbox{\raisebox{-0.01cm}{$\down{e}$}}\cdot u).\]
Hence the series expansion of a zonal measure $\mu$ becomes
\[\mu \sim \sum \limits_{k=0}^{\infty}N(n,k)\langle \mu,P_k^n(\mbox{\raisebox{-0.01cm}{$\down{e}$}}\cdot
.\:)\rangle P_k^n(\mbox{\raisebox{-0.01cm}{$\down{e}$}}\cdot
.\:).\]

\noindent The numbers $\mu_k:=\langle
\mu,P_k^n(\mbox{\raisebox{-0.01cm}{$\down{e}$}}\cdot .\:)\rangle$
are called Legendre coefficients of $\mu \in
\mathcal{M}(S^{n-1},\mbox{\raisebox{-0.01cm}{$\down{e}$}})$.
Using $\pi_k H=H$ for every $H \in \mathcal{H}_k^n$ and the fact,
that spherical convolution of zonal measures is commutative, we
obtain a version of the Funk-Hecke Theorem.

\begin{koro} If $\mu \in
\mathcal{M}(S^{n-1},\mbox{\raisebox{-0.01cm}{$\down{e}$}})$ and $H
\in \mathcal{H}_k^n$, then $H \ast \mu = \mu_k H$.
\end{koro}

We are now ready to give the definition of multiplier operators.

\begin{defi} \label{defmultiplier}
We call a map $\Phi: \mathcal{Q} \subseteq \mathcal{M}(S^{n-1})
\rightarrow \mathcal{M}(S^{n-1})$ a multiplier transformation if
there is a sequence of real numbers $c_k$ such that, for every $k
\in \mathbb{N}$,
\begin{equation} \label{defmulttrans}
\pi_k \Phi \mu = c_k \pi_k \mu, \qquad \forall \mu \in \mathcal{Q}
.
\end{equation}
The numbers $c_0, c_1, c_2, \ldots$ are called the multipliers of
$\Phi$.
\end{defi}

Using again the fact that spherical convolution of zonal measures
is commutative, we see that for $\mu \in
\mathcal{M}(S^{n-1},\mbox{\raisebox{-0.01cm}{$\down{e}$}})$ the
map $\Phi_{\mu}: \mathcal{M}(S^{n-1}) \rightarrow
\mathcal{M}(S^{n-1})$
\[\nu \mapsto \nu \ast \mu \]
is a multiplier transformation. The sequence of multipliers of
these convolution operators is just the sequence of Legendre
coefficients of the measure $\mu$.

By definition (\ref{defpik}) of the orthogonal projection $\pi_k$
and (\ref{rightconvintrot}), it is easy to see that multiplier
transformations intertwine rotations and that, by definition
(\ref{defmulttrans}), they are linear on the space
$\mathcal{H}^n$. The following corollary to Schur's Lemma
establishes the converse statement, see \cite{schneider74}, p.67.

\begin{satz} \label{multcrit} If $\Phi: \mathcal{H}^n \rightarrow
\mathcal{M}(S^{n-1})$ is an intertwining linear map, then $\Phi$
is a multiplier transformation.
\end{satz}

\vspace{1cm}

\setcounter{abschnitt}{3}
\centerline{\large{\bf{\arabic{abschnitt}. Convex Bodies and
Multiplier Transformations}}}

\reseteqn \alpheqn

\setcounter{koro}{0}

\vspace{0.7cm} \noindent We collect here further material on
convex geometry and endomorphisms of convex bodies, see
\cite{schneider93}, \cite{schneider74} and \cite{kiderlen05}. We
also prove that every Blaschke-Minkowski homomorphism is a
multiplier transformation.

The volume of a Minkowski linear combination $\lambda_1K_1 +
\ldots + \lambda_m K_m$ of convex bodies $K_1, \ldots, K_m$ is a
homogeneous polynomial of degree $n$ in the $\lambda_i$
\[V(\lambda_1K_1 + \ldots +\lambda_m K_m)=\sum \limits_{i_1,\ldots, i_n} V(K_{i_1},\ldots,K_{i_n})\lambda_{i_1}\cdots\lambda_{i_n}.   \]
The coefficients $V(K_{i_1},\ldots,K_{i_n})$ are called mixed
volumes of $K_{i_1}, \ldots, K_{i_n}$. These functionals are
nonnegative, translation invariant, monotone (with respect to set
inclusion) and multilinear with respect to Minkowski addition.
Denote by $V_i(K,L)$ the mixed volume $V(K,\ldots,K,L,\ldots,L)$,
where $K$ appears $n-i$ times and $L$ appears $i$ times. The
quermassintegrals $W_i(K)$ are given by $V_i(K,B^n)$.

Let $\mathcal{K}_i^n$ be the subset of $\mathcal{K}^n$ consisting
of convex bodies whose dimension is at least $n-i$. Then $K \in
\mathcal{K}_i^n$ if and only if $W_i(K) > 0$. The classical
inequality between two consecutive quermassintegrals states that
for $K \in \mathcal{K}^n$ and $0 \leq i \leq n-2$,
\begin{equation} \label{ineququermass}
W_{i+1}(K)^{n-i} \geq \kappa_n W_i(K)^{n-i-1},
\end{equation}
where $\kappa_n$ is the volume of the Euclidean unit ball $B^n$.
If $K \in \mathcal{K}_{i+1}^n$ there is equality in
(\ref{ineququermass}) if and only if $K$ is a ball.

\pagebreak

For the functional $V_1(K,L)$ there is an integral representation
\begin{equation} \label{v1intrep}
V_1(K,L)=\frac{1}{n}\langle h(L,\cdot),S_{n-1}(K,\cdot) \rangle.
\end{equation}
This shows that $V_1: [\mathcal{K}^n_0] \times \mathcal{K}^n
\rightarrow \mathbb{R}$ is bilinear with respect to Blaschke and
Minkowski addition.

Also the surface area measure of a Minkowski linear combination
of convex bodies $K_1, \ldots, K_m$ can be expressed as a
polynomial homogeneous of degree $n-1$
\begin{equation} \label{mixedsurfareameas}
S_{n-1}(\lambda_1K_1+\ldots +\lambda_mK_m,\cdot)=\sum
\limits_{i_1,\ldots,i_{n-1}} \lambda_{i_1}\cdots
\lambda_{i_{n-1}}S(K_{i_1},\ldots,K_{i_{n-1}},\cdot).
\end{equation}
The coefficients $S(K_{i_1},\ldots,K_{i_{n-1}},\cdot) \in
\mathcal{M}_o^+(S^{n-1})$ are called the mixed surface area
measures of $K_{i_1}, \ldots, K_{i_n}$. They have the property
that for any convex body $K$
\begin{equation} \label{mixedvolmixedsurf}
V(K,K_1,\ldots,K_{n-1})=\frac{1}{n} \langle
h(K,\cdot),S(K_1,\ldots,K_{n-1},\cdot) \rangle .
\end{equation}
The measures $S_j(K,\cdot):=S(K,\ldots,K,B^n,\ldots,B^n,\cdot)$,
where $K$ appears $j$ times and $B^n$ appears $n-1-j$ times, are
called the area measures of order $j$ of $K$.

By (\ref{pi0pi1}) and (\ref{mixedvolmixedsurf}), we have for $K
\in \mathcal{K}^n$,
\begin{equation} \label{wnminus1}
W_{n-1}(K)=\kappa_n \pi_0 h(K,\cdot) \qquad \mbox{and} \qquad
W_1(K)=\frac{1}{n}\pi_0 S_{n-1}(K,\cdot).
\end{equation}

A convex body $K \in \mathcal{K}^n$ is uniquely determined by its
support function $h(K,\cdot)$, which is positively homogeneous of
degree one and sublinear. Conversely, every function with these
properties is the support function of a convex body. By
(\ref{rotate1}), we have $\vartheta h(K,\cdot)=h(\vartheta
K,\cdot)$ for $\vartheta \in SO(n)$. Thus the support function of
a convex body $K$ is zonal if and only if $K$ is a body of
revolution. The Steiner point map $s: \mathcal{K}^n \rightarrow
\mathbb{R}^n$, defined by
\[s(K)=n\int_{S^{n-1}} h(K,u)u du,  \]
is up to a multiplicative constant the unique vector valued
continuous, rotation intertwining and Minkowski additive map, see
\cite{schneider74}. Since vector addition in $\mathbb{R}^n$
coincides with Minkowski addition of singletons, it is possible
to give an alternative definition of the Steiner point
\begin{equation} \label{Steinerpoint}
h(\{s(K)\},\cdot)=nh(K,\cdot) \ast
(\mbox{\raisebox{-0.01cm}{$\down{e}$}} \cdot .\:)=\pi_1h(K,\cdot).
\end{equation}

A convex body $K \in \mathcal{K}_0^n$ is also uniquely determined
up to translation by its surface area measure $S_{n-1}(K,\cdot)$
which is an element of $\mathcal{M}^+_o(S^{n-1})$, the set of
nonnegative measures on the sphere with center of mass in the
origin. Conversely, every element of $\mathcal{M}_o^+(S^{n-1})$
that is not concentrated on any great sphere is the surface area
measure of a convex body with interior points. For $\vartheta \in
SO(n)$, we have $\vartheta S_{n-1}(K,\cdot)=S_{n-1}(\vartheta
K,\cdot)$ and again the surface area measure of a convex body $K$
is zonal if and only if $K$ is a body of revolution.

There is no nonzero vector valued map from the set of translation
classes of convex bodies
$[\mathcal{K}_0^n]=\mathcal{K}_0^n/\mathbb{R}^n$ that is
continuous, rotation intertwining and additive with respect to
Blaschke addition. This fact is reflected by the relation
\begin{equation} \label{centerofmass}
\pi_1S_{n-1}(K,\cdot)=nS_{n-1}(K,\cdot) \ast
(\mbox{\raisebox{-0.01cm}{$\down{e}$}} \cdot .\:) =0.
\end{equation}

By (\ref{Steinerpoint}), the Steiner point map can be interpreted
as a convolution operator on the set of convex bodies. We will
consider in the following more general transformations induced by
convolution operators. By (\ref{convwithmeasure}), the convolution
from the left with measures $\mu \in \mathcal{M}^+(SO(n))$ can be
interpreted as (weighted) rotation means. The following
consequences of this interpretation appear in
\cite{grinbergzhang99}.

\begin{lem} \label{convvonlinks} Let $\mu \in \mathcal{M}^+(SO(n))$.
\begin{enumerate}
\item For $K \in \mathcal{K}^n$, the function $\mu \ast h(K,\cdot)$
is the support function of a convex body.
\item For $L \in \mathcal{K}_0^n$ and $\mu \neq 0$, the measure $\mu \ast
S_{n-1}(L,\cdot)$ is the surface area measure of a convex body
with interior points.
\end{enumerate}
\end{lem}

By (\ref{convright}) and the remarks after Definition
\ref{defmultiplier}, spherical convolution operators from the
right are multiplier transformations. It follows from
(\ref{centerofmass}) that the convolution of surface area measures
with a nonnegative zonal measure $\mu$ gives again nonnegative
measures with center of mass in the origin. It is also not hard
to see that $S_{n-1}(K,\cdot) \ast \mu$ is not concentrated on
any great sphere. Thus, the measure $S_{n-1}(K,\cdot) \ast \mu$ is
again a surface area measure of a convex body. Noting
(\ref{centerofmass}), we see that it is sufficient that the
measure $\mu$ is positive up to addition of a measure with
density $c (\mbox{\raisebox{-0.01cm}{$\down{e}$}} \cdot .\:)$. We
capture this property of a measure in the following definition:

\begin{defi} A measure $\mu \in \mathcal{M}(S^{n-1},\mbox{\raisebox{-0.01cm}{$\down{e}$}})$ is called weakly positive
if it is nonnegative up to addition of a linear measure, i.e.\ a
measure with density $c (\mbox{\raisebox{-0.01cm}{$\down{e}$}}
\cdot .\:)$, $c \in \mathbb{R}$.
\end{defi}

It was shown in \cite{kiderlen05} that also the cone of support
functions is invariant under convolution of zonal weakly positive
measures. We summarize these results in

\begin{lem} \label{supfaltsup} Let $\mu \in \mathcal{M}(S^{n-1},\mbox{\raisebox{-0.01cm}{$\down{e}$}})$ be weakly positive.
\begin{enumerate}
\item For $K \in \mathcal{K}^n$ the function $h(K,\cdot)
\ast \mu$ is the support function of a convex body.
\item For $L \in \mathcal{K}_0^n$ and $\mu$ not linear the measure
$S_{n-1}(L,\cdot) \ast \mu$ is the surface area measure of a
convex body with interior points.
\end{enumerate}
\end{lem}

We call a map $\Phi: \mathcal{K}^n \rightarrow \mathcal{K}^n$
that is continuous, rotation intertwining and Minkowski additive
a Minkowski endomorphism. A Blaschke endomorphism is a map $\Psi:
[\mathcal{K}_0^n] \rightarrow [\mathcal{K}_0^n]$ that is
continuous, rotation intertwining and additive with respect to
Blaschke addition.

Let $K, L \in \mathcal{K}^n$. Then $K \subseteq L$ if and only if
$h(K,\cdot) \leq h(L,\cdot)$. Thus a map $\Phi: \mathcal{K}^n
\rightarrow \mathcal{K}^n$ defined by
\[h(\Phi K,\cdot)=h(K,\cdot) \ast \mu,   \]
with a weakly positive measure $\mu \in
\mathcal{M}(S^{n-1},\mbox{\raisebox{-0.01cm}{$\down{e}$}})$ is, by
(\ref{Steinerpoint}), monotone (with respect to set inclusion) on
the set of convex bodies having their Steiner point in the
origin. We call a Minkowski endomorphism with this property
\emph{weakly monotone}.

\pagebreak

A classification of weakly monotone Minkowski endomorphisms and
Blaschke endomorphisms was established by Kiderlen in
\cite{kiderlen05}. We summarize his results in

\begin{satz} \label{MBendos} A map $\Phi: \mathcal{K}^n \rightarrow \mathcal{K}^n$ is a weakly monotone Minkowski
endomorphism if and only if there is a unique weakly positive
measure $\mu \in
\mathcal{M}(S^{n-1},\mbox{\raisebox{-0.01cm}{$\down{e}$}})$ such
that
\begin{equation} \label{minkendo}
h(\Phi K,\cdot)=h(K,\cdot) \ast \mu, \qquad K \in \mathcal{K}^n.
\end{equation}
A map $\Psi: [\mathcal{K}_0^n] \rightarrow [\mathcal{K}_0^n]$ is a
Blaschke endomorphism if and only if there is a weakly positive
measure $\nu \in
\mathcal{M}(S^{n-1},\mbox{\raisebox{-0.01cm}{$\down{e}$}})$,
unique up to addition of a linear measure, such that
\begin{equation} \label{blaschkeendo}
S_{n-1}(\Psi K,\cdot)= S_{n-1}(K,\cdot) \ast \nu, \qquad K \in
\mathcal{K}_0^n.
\end{equation}
\end{satz}

The major open problem concerning Minkowski endomorphisms is a
classification without the extra assumption of weak monotonicity.
For $n=2$, Schneider obtained in \cite{schneider74b} such a result
by showing that every Minkowski endomorphism is weakly monotone.
The following conjecture appears implicitly in \cite{schneider74b}
and \cite{kiderlen05}.

\begin{conj} \label{Mconj}
For $n \geq 3$ every Minkowski endomorphism is weakly monotone.
\end{conj}

In \cite{kiderlen05} a natural notion of adjointness between
Minkowski and Blaschke endomorphisms was introduced.

\begin{defi} A Minkowski endomorphism $\Phi$ and a Blaschke
endomorphism $\Psi$ are called adjoint if for every $K \in
\mathcal{K}_0^n$ and every $L \in \mathcal{K}^n$
\[V_1(\Psi K,L)=V_1(K,\Phi L).   \]
\end{defi}

Using (\ref{v1intrep}), Lemma \ref{selfadlemma} and Theorem
\ref{MBendos}, we see that a Blaschke and a Minkowski endomorphism
are adjoint if and only if they have the same generating measure
up to addition of a linear measure. By Theorem \ref{MBendos},
every Blaschke endomorphism has an adjoint weakly monotone
Minkowski endomorphism. The converse statement is equivalent to
Conjecture \ref{Mconj}.

The results obtained in Theorem \ref{MBendos} show that the
respective endomorphisms are multiplier transformations. This fact
has been deduced for Minkowski endomorphisms in \cite{schneider74}
using a different method. In the following we will adapt the
technique by Schneider to show that also Blaschke-Minkowski
homomorphisms induce multiplier transformations.

Every Blaschke-Minkowski homomorphism $\Phi: [\mathcal{K}_0^n]
\rightarrow \mathcal{K}^n$ induces a map on the set of surface
area measures by
\begin{equation} \label{induce}
\Phi S_{n-1}(K,\cdot)=h(\Phi K,\cdot), \qquad K \in
\mathcal{K}_0^n.
\end{equation}
Using Theorem \ref{multcrit}, we obtain:

\begin{satz} \label{conecond} Let $\Phi: [\mathcal{K}_0^n] \rightarrow \mathcal{K}^n$
be a Blaschke-Minkowski homomorphism. Then the induced map on the
set of surface area measures is a multiplier transformation, i.e.\
there is a sequence $c_k \in \mathbb{R}$ such that, for every $K
\in \mathcal{K}_0^n$,
\[\pi_k h(\Phi K,\cdot)=\pi_k \Phi S_{n-1}(K,\cdot)=c_k \pi_k S_{n-1}(K,\cdot).   \]
\end{satz}

For the proof of Theorem \ref{conecond}, we need some well known
facts on the vector space of differences of surface area measures,
see \cite{weil80} and \cite{groemer96}, p.70.

\begin{lem} \label{condsupport} Let $\mathcal{Q}\subseteq \mathcal{M}_o^+(S^{n-1})$ denote the set
of surface area measures of convex bodies with interior points.
Then
\begin{enumerate}
\item $\mathcal{Q}$ is dense in $\mathcal{M}_o^+(S^{n-1})$ and $\mathcal{M}_o(S^{n-1})=\mathcal{Q}-\mathcal{Q}.$
\item $\mathcal{Q} \cap \mathcal{H}^n$ is dense in $\mathcal{Q}$.
\end{enumerate}
\end{lem}

\noindent PROOF OF THEOREM \ref{conecond}: By the additivity
property of Blaschke-Minkowski homomorphisms, the induced map
(\ref{induce}) on the cone $\mathcal{Q}$ of surface area measures
of convex bodies is linear, and hence by Lemma \ref{condsupport}
(a), there is a unique linear extension $\tilde{\Phi}$ to the
vector space $\mathcal{M}(S^{n-1})$ given by
\[\tilde{\Phi}(\mu)=\Phi S_{n-1}(K_+,\cdot) - \Phi S_{n-1}(K_-,\cdot),\]
where $\mu-\pi_1 \mu=S_{n-1}(K_+,\cdot)-S_{n-1}(K_-,\cdot) \in
\mathcal{M}_o(S^{n-1})$ for some $K_+,K_- \in \mathcal{K}_0^n$.

The restriction of $\tilde{\Phi}$ to $\mathcal{H}^n$ is by
definition linear and intertwines rotations. Thus, by Theorem
\ref{multcrit}, it is a multiplier transformation. The result
follows since $\tilde{\Phi}$ and $\Phi$ coincide on the set
$\mathcal{Q} \cap \mathcal{H}^n$ which is dense in $\mathcal{Q}$
by Lemma \ref{condsupport} (b). \hfill $\blacksquare$\\

By Cauchy's surface area formula, the mean width of the projection
body of a convex body $K \in \mathcal{K}_0^n$ is a constant
multiple of the surface area of $K$. The following corollary to
Theorem \ref{conecond} is a generalization of this fact.

\begin{koro} \label{meanwidth} Let $\Phi: [\mathcal{K}_0^n] \rightarrow \mathcal{K}^n$
be a Blaschke-Minkowski homomorphism. Then
\[W_{n-1}(\Phi K)=r_{\Phi}W_{1}(K),  \]
where $r_{\Phi} \in \mathbb{R}^+$ is the radius of the ball $\Phi
B^n$.
\end{koro}
{\it Proof}: We will first show that $\Phi B^n$ is a ball. To see
this, note that $\pi_k S_{n-1}(B^n,\cdot)=0$ for $k \geq 1$. Thus
by Theorem \ref{conecond}, we have $\pi_k h(\Phi B^n,\cdot)=0$ for
$k \geq 1$, hence $\Phi B^n$ is a ball. By Theorem \ref{conecond},
the radius $r_{\Phi}$ of $\Phi B^n$ is given by
\[r_{\Phi}=\pi_0 h(\Phi B^n,\cdot)=\pi_0\Phi S_{n-1}(B^n,\cdot)=c_0\pi_0S_{n-1}(B^n,\cdot)=c_0 \omega_n,    \]
where $c_0$ denotes the first multiplier of $\Phi$ and $\omega_n$
is the surface area of $B^n$ . By (\ref{wnminus1}), we have
$W_{n-1}(\Phi K)=\kappa_n\pi_0 h(\Phi K)$ and thus, again by
Theorem \ref{conecond} and (\ref{wnminus1}),
\[W_{n-1}(\Phi
K) =\kappa_n\pi_0 \Phi S_{n-1}(K,\cdot) =\frac{r_{\Phi}}{n}\pi_0
S_{n-1}(K,\cdot)=r_{\Phi}W_1(K) .\]
\hfill $\blacksquare$ \\

\vspace{0.7cm}

\setcounter{abschnitt}{4}
\centerline{\large{\bf{\arabic{abschnitt}. Characterization of
Blaschke-Minkowski Homomorphisms}}}

\reseteqn \alpheqn

\setcounter{koro}{0}

\vspace{0.7cm} \noindent We turn now to the proofs of the main
theorems. From now on we will view  a map $\Phi: [\mathcal{K}_0^n]
\rightarrow \mathcal{K}^n$ via the obvious identification as a
translation invariant map on $\mathcal{K}_0^n$. The next lemma
shows that every Blaschke-Minkowski homomorphism has a unique
continuous extension to $\mathcal{K}^n$.

\begin{lem} \label{extension} Let $\Phi: \mathcal{K}_0^n \rightarrow \mathcal{K}^n$
be a Blaschke-Minkowski homomorphism. Then there is a unique
continuous extension of $\Phi$ to $\mathcal{K}^n$.
\end{lem}
{\it Proof}: Let $K_m \in \mathcal{K}_0^n$ be a sequence
converging to $K \in \mathcal{K}^n$. Then we define
\[\Phi K =\lim_{m \rightarrow \infty} \Phi K_m.    \]
To see that this limit exists, note that, by Corollary
\ref{meanwidth}, $W_{n-1}(\Phi K_m)=r_{\Phi}W_1(K_m)$. Thus, $
W_{n-1}(\Phi K_m) \rightarrow r_{\Phi}W_1(K)$ as $m \rightarrow
\infty$. Hence the sequence $\Phi K_m$ is bounded. Let $\Phi
K_{m_j}$ be a convergent subsequence of $\Phi K_m$ with limit $L
\in \mathcal{K}^n$. By Theorem \ref{conecond} and (\ref{defpik}),
\[\pi_k h(\Phi K_{m_j},\cdot) = c_k \pi_k S_{n-1}(K_{m_j},\cdot)=c_kN(n,k) S_{n-1}(K_{m_j},\cdot)
\ast P_k^n(\mbox{\raisebox{-0.01cm}{$\down{e}$}}  \cdot .\:).\]
By Lemma \ref{approxlem}, this converges uniformly to $c_k \pi_k
S_{n-1}(K,\cdot)$. On the other hand, $\pi_k h(\Phi K_{m_j},\cdot)
\rightarrow \pi_k h(L,\cdot)$ as $j \rightarrow \infty$. By the
completeness of spherical harmonics, this implies that the limit
of every convergent subsequence $\Phi K_{m_j}$ of $\Phi K_m$
coincides and thus $\Phi K_m$ itself is convergent.
\hfill $\blacksquare$ \\

We will need a criterion to determine if a measure $\mu \in
\mathcal{M}(S^{n-1},\mbox{\raisebox{-0.01cm}{$\down{e}$}})$ is
weakly positive. Let $\mathcal{L}=\{h(K,\cdot)-h(L,\cdot):K,L \in
\mathcal{K}^n\} $ denote the vector space of differences of
support functions. The following Lemma is in a slightly weaker
form due to Schneider \cite{schneider74b} for $n=2$ and Kiderlen
\cite{kiderlen05} for $n \geq 3$.

\begin{lem} \label{weakpos1} Let $g \in \mathcal{L}$ and let $\mathcal{N}$ be a
dense subset of $\mathcal{M}^+_o(S^{n-1})$. Then
\begin{equation} \label{gmu}
\langle g,\mu\rangle\geq 0 \qquad \forall\, \mu \in \mathcal{N}
\end{equation}
if and only if there is an $x \in \mathbb{R}^n$ such that
\begin{equation} \label{gweaklypos}
g(u)+x\cdot u \geq 0 \qquad \forall \, u \in S^{n-1}.
\end{equation}
\end{lem}
{\it Proof}: Obviously, (\ref{gweaklypos}) for some $x \in
\mathbb{R}^n$ implies (\ref{gmu}). Conversely, assume that
(\ref{gmu}) holds. Since $\mathcal{N}$ is dense in
$\mathcal{M}^+_o(S^{n-1})$, (\ref{gmu}) holds for every measure
in $\mathcal{M}^+_o(S^{n-1})$. Let
\[g=h(L,\cdot)-h(M,\cdot) \]
with convex bodies $L,M \in \mathcal{K}_0^n$. Define the inradius
of $L$ relative to $M$ by
\[r(L,M)=\max \{\lambda \geq 0: \lambda M \subseteq L + x \mbox{ for some } x \in \mathbb{R}^n\}. \]
Choose $x \in \mathbb{R}^n$, with $r(L,M)M \subseteq L + x$. By
the definition of $r(L,M)$, the contact points of $r(L,M)M$ and $L
+ x$ are distributed on their boundaries such that
\[o \in \mbox{conv}\{N(L,y) \cap S^{n-1}:y \in r(L,M)M \cap L+x   \},  \]
where $N(L,y)$ is the normal cone of $L$ in $y$. Otherwise we
could move the body $r(L,M)M$ inside $L + x$ away from the contact
points and blow it up, in contradiction to the definition of
$r(L,M)$. Let $\mu \in \mathcal{M}^+_o(S^{n-1})$ be concentrated
in the set $\{N(L,y) \cap S^{n-1}:y \in r(L,M)M \cap L+x\}$. By
(\ref{gmu}),
\[r(L,M)\langle h(M,\cdot),\mu \rangle=\langle h(L,\cdot),\mu \rangle=
\langle g+h(M,\cdot),\mu\rangle \geq \langle h(M,\cdot),\mu
\rangle. \] Thus $r(L,M) \geq 1$, and hence we have for every $u
\in S^{n-1}$
\[g(u)+h(M,u)+x\cdot u=h(L+x,u)\geq r(L,M)h(M,u) \geq h(M,u). \]
\hfill $\blacksquare$\\

Using (\ref{v1intrep}), and noting that the set of surface area
measures of convex bodies is a dense subset of
$\mathcal{M}^+_o(S^{n-1})$, we obtain the following geometric
consequence of Lemma \ref{weakpos1} which was proved differently
by Weil in \cite{weil74}.

\begin{koro} \label{weil} Let $K, L \in \mathcal{K}^n$. If $V_1(M,K) \leq V_1(M,L)$ for every $M \in \mathcal{K}_0^n$ then there is a vector $x \in
\mathbb{R}^n$ such that $K + x \subseteq L.$
\end{koro}

Note that, if the function $g$ in Lemma \ref{weakpos1} is zonal,
then the vector $x$ in (\ref{gweaklypos}) can be chosen as a
multiple of $\mbox{\raisebox{-0.01cm}{$\down{e}$}}$. The
following consequence of Lemma \ref{weakpos1}, which we will use
frequently, was also used in the proof of Theorem \ref{MBendos},
see \cite{kiderlen05}.

\begin{koro} \label{weakpos3} Let $\mu \in \mathcal{M}(S^{n-1},\mbox{\raisebox{-0.01cm}{$\down{e}$}})$ and let $\mathcal{N}$ be a
dense subset of $\mathcal{M}^+_o(S^{n-1})$. Then
\begin{equation} \label{numugeq0}
\nu \ast \mu \in \mathcal{M}^+_o(S^{n-1}) \qquad \forall \,\nu
\in \mathcal{N}
\end{equation}
if and only if $\mu$ is weakly positive.
\end{koro}
{\it Proof}: It is clear that (\ref{numugeq0}) holds if $\mu$ is
weakly positive. Conversely, assume that (\ref{numugeq0}) holds.
Let $(\phi_k)_{k \in \mathbb{N}}$ be a zonal approximate identity.
Then $\nu \ast \mu \ast \phi_k \geq 0$, and by Lemma
\ref{zonalapproxid}, $\mu \ast \phi_k \in
\mathcal{C}^{\infty}(S^{n-1}).$ Using (\ref{zonalconv}), we see
that $(\nu \ast \mu \ast
\phi_k)(\mbox{\raisebox{-0.01cm}{$\down{e}$}})=\langle \mu \ast
\phi_k,\nu \rangle \geq 0$ for every $\nu \in \mathcal{N}$. As
$\mathcal{C}^{\infty}(S^{n-1})\subseteq \mathcal{L}$, see
\cite{schneider93} p.27, by Lemma \ref{weakpos1} and the remark
after Corollary \ref{weil}, there are $c_k \in \mathbb{R}$ such
that
\[(\mu \ast \phi_k)(u)+c_k(\mbox{\raisebox{-0.01cm}{$\down{e}$}}\cdot u) \geq 0.   \]
Thus, for nonnegative $f \in \mathcal{C}(S^{n-1})$, we have by
(\ref{pi0pi1})
\[f \ast \mu \ast \phi_k \geq -c_k f \ast (\mbox{\raisebox{-0.01cm}{$\down{e}$}}\cdot .)=-\frac{c_k}{n}\pi_1 f.   \]
By Lemma \ref{zonalapproxid}, $\mu \ast \phi_k \rightarrow \mu$
weakly, and thus $f \ast \mu \ast \phi_k \rightarrow f \ast \mu$
uniformly by Lemma \ref{approxlem}. Hence there exists $b \in
\mathbb{R}$ such that $b \geq -c_k \pi_1f$. Since $\pi_1f$ is a
linear \linebreak functional, the sequence $c_k$ is bounded.
Therefore we can assume that $c_k
\rightarrow c$. \hfill $\blacksquare$\\

The main ingredient in the proof of Theorem \ref{satzbmhomo} is a
classification of translation invariant homogeneous valuations of
convex sets. Since the map $K \mapsto S_{n-1}(K,\cdot)$ is a
translation invariant valuation, see \cite{schneider93}, p.201,
we obtain from the definition of Blaschke addition that for all
$K, L \in \mathcal{K}^n_0$ such that $K \cup L \in
\mathcal{K}_0^n$ and $K \cap L \in \mathcal{K}_0^n$,
\begin{equation} \label{blaschkerel}
(K \cup L) \: \# \: (K \cap L)=K \: \# \: L.
\end{equation}
Thus, if $\Phi$ is a Blaschke-Minkowski homomorphism, we have by
Lemma \ref{extension} for all $K, L \in \mathcal{K}^n$ such that
$K \cup L \in \mathcal{K}^n$,
\begin{equation} \label{bmminkval}
\Phi(K \cup L) + \Phi(K \cap L)=\Phi K + \Phi L.
\end{equation}
Hence, $\Phi$ is a valuation with respect to Minkowski addition.
For further information on valuations of this type see
\cite{ludwig02} and \cite{ludwig04}.

The following characterization is due to Hadwiger
\cite{hadwiger51} and McMullen \cite{mcmullen80}:

\begin{satz} \label{mcmullen} A map $\phi: \mathcal{K}^n \rightarrow \mathbb{R}$ is a continuous translation invariant
valuation homogeneous of degree $n-1$ if and only if there is a
function $g \in \mathcal{C}(S^{n-1})$, unique up to addition of a
linear function, such that
\[\phi(K)=\langle g,S_{n-1}(K,\cdot)\rangle.   \]
\end{satz}

Using Theorem \ref{mcmullen} and (\ref{bmminkval}), we can derive
the representation theorem for Blaschke-Minkowski homomorphisms.

\vspace{0.5cm} \noindent PROOF OF THEOREM \ref{satzbmhomo}:
Define a functional $\phi: \mathcal{K}^n \rightarrow \mathbb{R}$
by
\[\phi(K)=h(\Phi K,\mbox{\raisebox{-0.01cm}{$\down{e}$}}).   \]
Since $S_{n-1}(\lambda K,\cdot)=\lambda^{n-1}S_{n-1}(K,\cdot)$ for
$\lambda \geq 0$ and $K \in \mathcal{K}^n$, we have by
(\ref{bmadd})
\begin{equation} \label{bmhomogeneous}
\Phi \lambda K = \lambda^{n-1} \Phi K.
\end{equation}
Using (\ref{bmhomogeneous}) and (\ref{bmminkval}), we see that the
map $\phi$ is a continuous valuation on $\mathcal{K}^n$
homogeneous of degree $n-1$. By Theorem \ref{mcmullen}, there is a
function $g \in \mathcal{C}(S^{n-1})$, unique up to addition of a
linear function, such that
\[\phi(K)=\langle g,S_{n-1}(K,\cdot)\rangle.   \]
Since $\phi$ is invariant under rotations leaving
$\mbox{\raisebox{-0.01cm}{$\down{e}$}}$ fixed, the function $g$ is
zonal, and thus, by (\ref{rotate1}) and (\ref{rotate2}),
\begin{equation} \label{drehmittel}
h(\Phi K,\mbox{\raisebox{-0.06cm}{$\down{\eta}$}})=h(\Phi K,\eta
\mbox{\raisebox{-0.01cm}{$\down{e}$}} )=\langle
g,S_{n-1}(\eta^{-1}K,\cdot)\rangle =\langle \eta
g,S_{n-1}(K,\cdot) \rangle.
\end{equation}
(\ref{bmhomorep}) follows now from (\ref{zonalconv}) and
(\ref{zonalconv2}). To see that $g$ is weakly positive, note that
by (\ref{Steinerpoint}), (\ref{centerofmass}) and the
commutativity of the convolution of zonal functions,
\[h(\{s(\Phi K)\},\cdot)=nh(\Phi K,\cdot) \ast (\mbox{\raisebox{-0.01cm}{$\down{e}$}} \cdot .\:)
=nS_{n-1}(K,\cdot) \ast (\mbox{\raisebox{-0.01cm}{$\down{e}$}}
\cdot .\:) \ast g = 0. \] Since $s(\Phi K) \in \mbox{relint}\,
\Phi K$, see \cite{schneider93}, p.43,  we have $h(\Phi
K,\cdot)\geq 0$. Thus, noting that the set of surface area
measures is a dense subset of $\mathcal{M}_o^+(S^{n-1})$, it
follows from Corollary \ref{weakpos3}
that $g$ is weakly positive. \hfill $\blacksquare$\\

For later applications, we state further properties of the
generating functions of Blaschke-Minkowski homomorphisms in the
following Lemma.

\begin{lem} \label{propg} Let $g \in \mathcal{C}(S^{n-1},\mbox{\raisebox{-0.01cm}{$\down{e}$}})$ be the generating
function of a Blaschke-Minkowski homomorphism.
\begin{enumerate}
\item $g$ is a difference of support functions, i.e.\ $g \in \mathcal{L}.$
\item There is a symmetric body of revolution $L \in
\mathcal{K}^n$, such that for every $u \in S^{n-1}$,
\[g(u)+g(-u)=h(L,u).  \]
\end{enumerate}
\end{lem}
{\it Proof}: By Lemma \ref{condsupport} (a), there are convex
bodies $K_+, K_- \in \mathcal{K}_0^n$ such that
\[\delta_{\down{e}} - \pi_1\delta_{\down{e}} =S_{n-1}(K_+,\cdot)-S_{n-1}(K_-,\cdot).   \]
Since the Dirac measure $\delta_{\down{e}}$ is the neutral element
for zonal convolution, and as $(\pi_1 \delta_{\down{e}})(u)=n
\mbox{\raisebox{-0.01cm}{$\down{e}$}} \cdot u$ by (\ref{pi0pi1}),
we obtain
\[(\delta_{\down{e}}-\pi_1\delta_{\down{e}})\ast g = g-\pi_1 g=h(\Phi K_+,\cdot)-h(\Phi K_-,\cdot).   \]
Since $\pi_1 g$ is a linear functional on $\mathbb{R}^n$, there is
a vector $x \in \mathbb{R}^n$ such that
\[(\pi_1g)(u)=x \cdot u=h(\{x\},u).  \]
Hence $g=h(\Phi K_++x,\cdot)-h(\Phi K_-,\cdot)$  which proves (a).

To see (b), let $\{b_1, \ldots, b_n\}$ be an orthonormal basis in
$\mathbb{R}^n$ such that
$\mbox{\raisebox{-0.01cm}{$\down{e}$}}=b_n$. For a vector $x \in
\mathbb{R}^n$ let $x_1,\ldots,x_n$ denote its coordinates with
respect to $b_1,\ldots,b_n$. Choose $\beta \in \mathbb{R}^+$ such
that the ellipsoid $E_{\alpha}$ defined by
\[\frac{x_1^2+\ldots+x_{n-1}^2}{\alpha^2}+\frac{x_n^2}{\beta^2}\leq 1  \]
has surface area $S(E_{\alpha})=1$. It was shown in
\cite{grinbergzhang99}, p.103, that as $\alpha \rightarrow
\infty$, we have $\beta \rightarrow 0$ and
\[S_{n-1}(E_{\alpha},\cdot) \rightarrow \frac{1}{2}(\delta_{\down{e}}+\delta_{-\down{e}})   \]
weakly. By Lemma \ref{approxlem},
\[h(\Phi E_{\alpha},u)=(S_{n-1}(E_{\alpha},\cdot) \ast g)(u) \rightarrow \frac{1}{2}(g(u)+g(-u))   \]
uniformly in $u \in S^{n-1}$. Since $h(\Phi E_{\alpha},\cdot)$
converges uniformly, it converges to
a support function of a convex body, which proves (b). \hfill $\blacksquare$\\

An immediate consequence of Lemma \ref{propg} is the complete
classification of all even Blaschke-Minkowski homomorphisms.

\vspace{0.5cm} \noindent PROOF OF THEOREM \ref{bmsymm}: A
Blaschke-Minkowski homomorphism is even if and only if its
generating function is even. Thus the result
follows from Lemma \ref{propg} (b). \hfill $\blacksquare$\\

If $g=h(L,\cdot)$ for some body of revolution $L \in
\mathcal{K}^n$ is the generating function of a Blaschke-Minkowski
homomorphism $\Phi$, then by (\ref{drehmittel}) and
(\ref{v1intrep}),
\begin{equation} \label{mixedvolint}
h(\Phi K,\mbox{\raisebox{-0.06cm}{$\down{\eta}$}})=nV_1(K,\eta L).
\end{equation}
Since $K_1 \subseteq K_2$ if and only if $h(K_1,\cdot) \leq
h(K_2,\cdot)$, the monotonicity of mixed volumes together with
(\ref{mixedvolint}) implies

\begin{koro} \label{bmmonotone} A Blaschke-Minkowski homomorphism
whose generating function is given by $h(L,\cdot)$ for some $L
\in \mathcal{K}^n$, is monotone with respect to set inclusion.
\end{koro}

Note that by Theorem \ref{bmsymm} and Corollary \ref{bmmonotone},
every even Blaschke-Minkowski homomorphism is monotone.

By Lemma \ref{convvonlinks}, every map $\Phi: \mathcal{K}^n
\rightarrow \mathcal{K}^n$ of the form
\[h(\Phi K,\cdot)=S_{n-1}(K,\cdot) \ast h(L,\cdot)   \]
for some $L \in \mathcal{K}^n$ is  a Blaschke-Minkowski
homomorphism, but in general there are generating functions $g$
of Blaschke-Minkowski homomorphisms that are not support
functions. An example of such a map is the (normalized) second
mean section operator $M_2$ introduced in \cite{goodeyweil92} and
further investigated in \cite{schneiderhug02}: Let
$\mathcal{E}_2^n$ be the affine Grassmanian of two-dimensional
planes in $\mathbb{R}^n$ and $\mu_2$ its motion invariant
measure, normalized such that $\mu_2(\{E \in \mathcal{E}_2^n: E
\cap B^n \neq \varnothing\})=\kappa_{n-2}$. Then
\[h(M_2K,\cdot)=(n-1)\int \limits_{\mathcal{E}_2^n} h(K \cap E,\cdot)d\mu_2(E)-h(\{z_{n-1}(K)\},\cdot)=S_{n-1}(K,\cdot) \ast g_2,   \]
where $z_{n-1}(K)$ is the intrinsic $(n-1)$st moment vector of
$K$, see \cite{schneider93}, p.304, and $\Lambda g_2$ is given by
\[\Lambda g_2(t)=\arccos(-t)\sqrt{1-t^2}.  \]
The function $g_2$ is not a support function. Note that the
operator $M_2$ is not monotone but has the following weak
monotonicity property: $M_2$ is monotone on those convex bodies
having their $(n-1)$st intrinsic moment vector in the origin.
This is similar to the monotonicity property of weakly monotone
Minkowski endomorphisms.

We will give now a complete characterization of generating
functions of Blaschke-Minkowski homomorphisms in the spirit of a
classification result of Weil \cite{weil76} of generating
measures of generalized zonoids. To this end, we need the
extension of area measures of convex bodies to the space
$\mathcal{L}$ of differences of support functions.

\begin{defi} Let $g_i \in \mathcal{L}, \: i=1,\ldots, n-1$, with
$g_i=h(K_i^0,\cdot)-h(K_i^1,\cdot)$. Then the mixed surface area
measure of $g_1,\ldots, g_{n-1}$ is defined by
\[S(g_1,\ldots,g_{n-1},\cdot)=\sum \limits_{\alpha_1,\ldots,\alpha_{n-1} \in \{0, 1\}}
 (-1)^{\alpha_1+\ldots+\alpha_{n-1}}S(K_1^{\alpha_1},\ldots,K_{n-1}^{\alpha_{n-1}},\cdot) \in \mathcal{M}_o(S^{n-1}).   \]
For a function $f \in \mathcal{C}(S^{n-1})$, define
\[V(f,g_1,\ldots,g_{n-1})=\langle f,S(g_1,\ldots,g_{n-1},\cdot)\rangle.   \]
For $g \in \mathcal{L}$ and $j=1,\ldots, n-1$, the measure
$S_j(g,\cdot)=S(g,\ldots,g,1,...,1,\cdot)$, where $g$ appears $j$
times and $1$ appears $n-j-1$ times, is called the area measure
of order $j$ of $g$.
\end{defi}

If $\Phi$ is a Blaschke-Minkowski homomorphism, then by Lemma
\ref{propg} (a),
\[ h(\Phi K,\cdot)=S_{n-1}(K,\cdot) \ast g = S_{n-1}(K,\cdot)
\ast h(L_+,\cdot) - S_{n-1}(K,\cdot) \ast h(L_-,\cdot),\] where
$g=h(L_+,\cdot)-h(L_-,\cdot)$. Thus, defining Blaschke-Minkowski
homomorphisms $\Phi_+$ and $\Phi_-$ with generating functions
$h(L_+,\cdot)$ and $h(L_-,\cdot)$ we obtain
\begin{equation} \label{diffsup}
h(\Phi K,\cdot) = h(\Phi_+K,\cdot)-h(\Phi_-K,\cdot).
\end{equation}
In the light of (\ref{diffsup}), we need a criterion to determine
whether a difference of support functions is in fact a support
function. This was established by Weil in \cite{weil74b}.

\begin{satz} \label{diffsupistsup} A function $g \in \mathcal{L}$ is the support
function of a convex body $K$ if and only if for all $j \in
\{1,\ldots,n-1\}$,
\[S_j(g,\cdot) \in \mathcal{M}_o^+(S^{n-1}).   \]
\end{satz}

In order to use Theorem \ref{diffsupistsup}, we need to determine
the area measures $S_j(\Phi K,\cdot)$. In \cite{grinbergzhang99},
p.105, the area measures of the convex body with support function
$\mu \ast h(K,\cdot)$, $\mu \in \mathcal{M}^+(SO(n))$ were
calculated. The result established there extends easily to
differences of support functions. Identifying spherical measures
with right $SO(n-1)$ invariant measures on $SO(n)$, we get the
following Lemma.

\begin{lem} \label{areameasurebm} Let $\Phi$ be a Blaschke-Minkowski homomorphism with
generating function $g \in \mathcal{L}$. Then $\langle f,
S_j(\Phi K,\cdot) \rangle$ is given by
\[\int_{(S^{n-1})^j} V(f,\Lambda g(u_1 \cdot .),\ldots,\Lambda g(u_j\cdot
.),1,\ldots,1)dS_{n-1}(K,u_1)\ldots dS_{n-1}(K,u_j).\]
\end{lem}

Using Lemma \ref{areameasurebm}, Theorem \ref{diffsupistsup} and
the fact that the set of surface area measures of convex bodies
forms a dense subset of $\mathcal{M}_o^+(S^{n-1})$, we obtain the
following characterization of generating functions of
Blaschke-Minkowski homomorphisms.

\begin{satz} A function $g \in \mathcal{L}$ is the generating function of a
Blaschke-Minkowski homomorphism if and only if, for every $j
=1,\ldots,n-1$,
\[\int_{(S^{n-1})^j} V(f,\Lambda g(u_1 \cdot .),\ldots,\Lambda g(u_j\cdot
.),1,\ldots,1)d\mu(u_1)\ldots d\mu(u_j)\geq 0\] for every
nonnegative $f \in \mathcal{C}(S^{n-1})$ and every $\mu \in
\mathcal{M}^+_o(S^{n-1})$.
\end{satz}

\vspace{1cm}

\setcounter{abschnitt}{5}
\centerline{\large{\bf{\arabic{abschnitt}. Endomorphisms and
Homomorphisms of Convex Bodies}}}

\reseteqn \alpheqn

\setcounter{koro}{0}

\vspace{0.7cm} \noindent We turn now to the connection between
adjoint Minkowski and Blaschke endomorphisms and
Blaschke-Minkowski homomorphisms.

\vspace{0.5cm} \noindent PROOF OF THEOREM \ref{bmanddual}: If
$\Psi$ and $\Psi^*$ are adjoint, then $\Psi$ is weakly monotone
and they have the same generating measure $\mu \in
\mathcal{M}(S^{n-1},\mbox{\raisebox{-0.01cm}{$\down{e}$}})$. Let
$\Phi$ be a Blaschke-Minkowski homomorphism with generating
function $g \in
\mathcal{C}(S^{n-1},\mbox{\raisebox{-0.01cm}{$\down{e}$}})$. From
the commutativity of zonal convolution, it follows that
\begin{eqnarray*}
h(\Phi \Psi^*K,\cdot)& = & S_{n-1}(\Psi^*K,\cdot) \ast
g=S_{n-1}(K,\cdot) \ast \mu \ast g \\
& = & S_{n-1}(K,\cdot) \ast g \ast \mu = h(\Phi K,\cdot) \ast \mu
= h(\Psi \Phi K,\cdot).
\end{eqnarray*}
Thus (a) implies (b) and obviously (b) implies (c).

By the multiplier property, a Blaschke-Minkowski homomorphism
$\Phi$ is injective if and only if all the multipliers of $g$ are
nonzero. Thus, the multipliers of $\Psi^*$ and $\Psi$ can be
determined from $\Phi \circ \Psi^*$ and $\Psi \circ \Phi$ and are
equal if (\ref{phipsipsiphi}) holds. By the completeness of the
system of spherical harmonics, it follows that (c) implies (a). \hfill $\blacksquare$\\

\noindent Theorem \ref{bmanddual} shows that the following
conjecture is equivalent to Conjecture \ref{Mconj}:

\begin{conj}
There exists an injective Blaschke-Minkowski homomorphism whose
range is invariant under every Minkowski endomorphism.
\end{conj}

\vspace{0.3cm} In view of this formulation of Conjecture
\ref{Mconj}, we further investigate the range of
Blaschke-Minkowski homomorphisms.

\begin{satz} The range of every Blaschke-Minkowski homomorphism is nowhere dense in $\mathcal{K}^n$.
\end{satz}
{\it Proof}: We call $K \in \mathcal{K}_0^n$ Blaschke
decomposable if there exist two bodies $K_1, K_2 \in
\mathcal{K}_0^n$ not homothetic to $K$ such that $K=K_1 \: \# \:
K_2$. By a result of Bronshtein \cite{bronshtein79}, the only
Blaschke indecomposable bodies in $\mathcal{K}_0^n$ are the
simplices. Thus, every body in the range of a Blaschke-Minkowski
homomorphism with the only possible exception of the image of
simplices is decomposable with respect to Minkowski addition.

Since the image of simplices is nowhere dense in $\mathcal{K}^n$
and since, on the other hand, the indecomposable bodies with
respect to Minkowski addition form a dense subset of
$\mathcal{K}^n$, the desired result
follows. \hfill $\blacksquare$\\

In the second part of this section we will see that most of the
geometric convolution operators we encountered so far do not
attain values in the set of polytopes.

\begin{satz} \label{bmsuppnopoly2} Let $\Phi: \mathcal{K}^n \rightarrow \mathcal{K}^n$
be a Blaschke-Minkowski homomorphism generated by the support
function $h(L,\cdot)$ of a body of revolution $L \in
\mathcal{K}^n$. If there is a convex body $K \in \mathcal{K}_0^n$
such that $\Phi K$ is a polytope, then there is a constant $c \in
\mathbb{R}^+$ such that
\[\Phi = c \Pi.   \]
\end{satz}
{\it Proof}: Let $P=\Phi K=\mbox{conv}\{x_1,\ldots,x_k\}$ be a
polytope with vertices $x_1, \ldots, x_k$. Then
\[h(P,\cdot)=S_{n-1}(K,\cdot) \ast h(L,\cdot).   \]
Since the body $L \in \mathcal{K}^n$ is unique up to translation,
we can assume that $h(L,\cdot) \geq 0$. Let $\mu =
\up{S}_{n-1}(K,\cdot) \in \mathcal{M}^+(SO(n))$, then by
(\ref{bmhomorep})
\begin{equation} \label{rotmean}
h(P,\cdot)=\int_{SO(n)} h(\vartheta L,\cdot) d\mu(\vartheta).
\end{equation}
From now on, we consider support functions as positive homogeneous
functions on $\mathbb{R}^n$. Let $C_1, \ldots, C_k$ denote the
normal cones of the vertices of $P$. Then the support function
$h(P,\cdot)$ is linear in every $C_i, \: i=1,\ldots,k$. Thus, by
(\ref{rotmean}), we have
\begin{equation} \label{gag}
\int_{SO(n)} h(\vartheta L,v_1)+h(\vartheta L,v_2)-h(\vartheta
L,v_1+v_2)d\mu(\vartheta)=0
\end{equation}
for all $v_1, v_2 \in C_i$. Since support functions are sublinear,
the integrand in (\ref{gag}) is nonnegative. Thus, as $\mu$ is
nonnegative, $h(\vartheta L,v_1)+h(\vartheta L,v_2)=h(\vartheta
L,v_1+v_2)$  for all $\vartheta$ in the support of $\mu$. For
each such $\vartheta$, we thus have
\[h(L,v_1)+h(L,v_2)=h(L,v_1+v_2)  \]
for all $v_1,v_2 \in \vartheta C_i$. This implies that $L$ is a
polytope itself. But since $L$ is a body of revolution and the
only polytopes that are bodies of revolution are the multiples of
the segment
$[-\mbox{\raisebox{-0.01cm}{$\down{e}$}},\mbox{\raisebox{-0.01cm}{$\down{e}$}}]$,
the desired result
follows from (\ref{projbody}). \hfill $\blacksquare$\\

\noindent Note that Theorem \ref{bmsymm} and Theorem
\ref{bmsuppnopoly2} imply Theorem \ref{bmsuppnopoly}. For a
corresponding result in dimension two see \cite{schneider74b},
p.311.

The Difference body operator $D: \mathcal{K}^n \rightarrow
\mathcal{K}^n$ is the Minkowski endomorphism defined by
\[DK=K+(-K).   \]
The Blaschke body operator $\nabla: \mathcal{K}^n \rightarrow
\mathcal{K}^n$ is the Blaschke endomorphism defined by
\[\nabla K=K \: \# \:(-K).   \]

\begin{koro} The only even Blaschke endomorphisms taking values in the set
of polytopes are constant multiples of $\nabla$.

If an even Minkowski endomorphism maps a zonoid onto a polytope,
then it is a constant multiple of $D$.
\end{koro}
{\it Proof}: Let $\Psi$ be an even Blaschke endomorphism and let
$\Psi K=P$ be a polytope for some $K \in \mathcal{K}_0^n$. By
(\ref{bmadd}), the map $\Pi \circ \Psi$ is an even
Blaschke-Minkowski homomorphism such that $\Pi \Psi K$ is a
polytope. By Theorem \ref{bmsuppnopoly} and Theorem \ref{bmsymm},
there is a constant $c \in \mathbb{R}^+$ such that
\begin{equation} \label{pipsipi}
\Pi \circ \Psi = c\Pi.
\end{equation}
Since $\Pi$ is injective all the even multipliers of $\Pi$ are
nonzero. Thus, by (\ref{pipsipi}), all even multipliers of $\Psi$
are equal to $c$. Noting that the odd multipliers of even
multiplier operators are zero, the result follows.

An analogous argument leads to the second statement. \hfill
$\blacksquare$

\vspace{1cm}

\setcounter{abschnitt}{6}
\centerline{\large{\bf{\arabic{abschnitt}. Geometric Inequalities
and Induced Operators}}}

\reseteqn \alpheqn

\setcounter{koro}{0}

\vspace{0.7cm} \noindent An important open problem in the theory
of affine isoperimetric inequalities is the conjectured
projection inequality by Petty \cite{petty72}:
\begin{equation} \label{pettyconj}
\frac{\kappa_{n}^{n-1}}{\kappa_{n-1}^n}V(\Pi K) \geq \kappa_n
V(K)^{n-1}
\end{equation}
with equality if and only if $K$ is an ellipsoid. If
(\ref{pettyconj}) holds, then, as was shown in \cite{lutwak90b},
it is a strengthened version of the classical isoperimetric
inequality $W_1(K)^n \geq \kappa_n V(K)^{n-1}$, compare Lemma
\ref{upbound} and (\ref{schconj}).

In this chapter we will study analogous problems for general
Blaschke-Minkowski homomorphisms and related operators which will
be introduced in the next theorem. Most of the results in this
chapter were established for the projection body operator in
\cite{lutwak90b}, see also \cite{lutwak85}, \cite{lutwak93}. The
aim of this section is to generalize the results obtained there
to general (nontrivial) Blaschke-Minkowski homomorphisms and to
show that the crucial tool is a representation of the form of
(\ref{bmhomorep}).

In the following $\Phi: \mathcal{K}^n \rightarrow \mathcal{K}^n$
shall always denote a {\it nontrivial} Blaschke-Minkowski
homomorphism.

\begin{satz} The map $\Phi$ satisfies the Steiner type formula
\[\Phi(K+\varepsilon B^n)=\sum \limits_{i=0}^{n-1}\varepsilon^{i} {n-1 \choose i}\Phi_{i}K.   \]
The operators $\Phi_i: \mathcal{K}^n \rightarrow \mathcal{K}^n,
\: i=0,\ldots, n-1,$ are continuous, translation invariant,
rotation intertwining Minkowski valuations and $\Phi_0=\Phi$.
\end{satz}
{\it Proof}: The desired result is an immediate consequence of
Theorem \ref{satzbmhomo} and the Steiner formula for the surface
area measure of a convex body $K$, see (\ref{mixedsurfareameas}),
\[S_{n-1}(K+\varepsilon B^n,\cdot)=\sum \limits_{i=0}^{n-1} \varepsilon^i {n-1 \choose i} S_{n-1-i}(K,\cdot) .\]
If $g \in
\mathcal{C}(S^{n-1},\mbox{\raisebox{-0.01cm}{$\down{e}$}})$
denotes the generating function of $\Phi$, then
\begin{equation} \label{defphii}
h(\Phi_i K,\cdot)=S_{n-1-i}(K,\cdot) \ast g.
\end{equation}
By Minkowski's existence theorem, the area measure $S_i(K,\cdot)$
of order $i$ of a convex body $K \in \mathcal{K}_0^n$ is also the
surface area measure $S_{n-1}(L,\cdot)$ of order $n-1$ of some
convex body $L \in \mathcal{K}_0^n$. Thus, the mappings $\Phi_i$
are well defined. Since the mappings $K \mapsto S_i(K,\cdot)$ are
translation invariant valuations, the operators $\Phi_i$ are
Minkowski valuations. \hfill $\blacksquare$ \\

Note that, by Theorem \ref{satzbmhomo} and
(\ref{mixedsurfareameas}), the mappings $\Phi_i$ are special cases
of more general operators defined on the cartesian product of
$n-1$ copies of $\mathcal{K}^n$. These mappings are studied in
more detail in \cite{schuster05}. In the following we will
consider only the operators $\Phi_i, \: i=0, \ldots, n-2,$ since
$\Phi_{n-1}$ maps every body $K$ to $\Phi B^n$ because
$S_0(K,\cdot)=S_{n-1}(B^n,\cdot)$ is independent of $K$. We
remark here that, for $K \in \mathcal{K}^n_i$, the image $\Phi_i
K$ is an element of $\mathcal{K}_0^n$ and $\Phi_i L=o$ if $L \in
\mathcal{K}^n_{i+2}\backslash \mathcal{K}^n_{i+1}$.

By (\ref{defphii}), the $\Phi_i$ are multiplier operators, but
apart from $\Phi_{0}=\Phi$ and $\Phi_{n-2}$ they can not be
interpreted as additive transformations of convex bodies, since
the set of area measures $S_j(K,\cdot)$ of order $j$  does not
form a cone in $\mathcal{M}^+_o(S^{n-1})$ for $j=2,\ldots,n-2$,
see \cite{goodeyschneider80}. The operator $\Phi_{n-2}$ is a
Minkowski endomorphism. To see this, note that the area measure
$S_1(K,\cdot)$ of  order one is related to the support function
$h(K,\cdot)$ by the linear second order differential operator
\[\Delta_1= \Delta_0 + (n-1),    \]
where $\Delta_0$ denotes the Laplace Beltrami operator on
$S^{n-1}$, see \cite{grinbergzhang99}, p.87. We have
\begin{equation} \label{delta1hk}
\Delta_1 h(K,\cdot)=S_1(K,\cdot),
\end{equation}
where this equality is understood in the sense of distributions
if $h(K,\cdot)$ is not in $\mathcal{C}^2(S^{n-1})$. From
(\ref{delta1hk}), it follows that
$S_1(K_1+K_2,\cdot)=S_1(K_1,\cdot)+S_1(K_2,\cdot)$, which together
with (\ref{defphii}) shows that $\Phi_{n-2}$ is a Minkowski
endomorphism.

As $\Delta_0$ is an intertwining operator so is $\Delta_1$. Thus,
by Lemma \ref{multcrit}, $\Delta_1$ is a multiplier operator. For
the following Lemma see \cite{grinbergzhang99}, p.86, and note
that multiplier transformations are obviously commutative.

\begin{lem} \label{delta1lem} Let $\mu \in
\mathcal{M}(S^{n-1},\mbox{\raisebox{-0.01cm}{$\down{e}$}})$ and
$\nu \in \mathcal{M}(S^{n-1})$. Then
\[\Delta_1(\nu \ast \mu)=\nu \ast (\Delta_1 \mu) = (\Delta_1 \nu) \ast \mu   \]
in the sense of distributions.
\end{lem}

\noindent Using Lemma \ref{delta1lem}, we get the following
result.

\begin{satz} The operator
$\Phi_{n-2}$ is a weakly monotone Minkowski endomorphism.
\end{satz}
{\it Proof}: We have seen that $\Phi_{n-2}$ is a Minkowski
endomorphism. In order to prove that $\Phi_{n-2}$ is weakly
monotone, we need to show, by Theorem \ref{MBendos}, that there is
a weakly positive measure $\mu \in
\mathcal{M}(S^{n-1},\mbox{\raisebox{-0.01cm}{$\down{e}$}})$ such
that $h(\Phi_{n-2}K,\cdot)=h(K,\cdot) \ast \mu.$

If $g \in \mathcal{C}(S^{n-1},\raisebox{-0.01cm}{$\down{e}$})$ is
the generating function of $\Phi$, then by Lemma \ref{delta1lem}
and (\ref{delta1hk}),
\[h(\Phi_{n-2}K,\cdot)=S_1(K,\cdot) \ast g = h(K,\cdot) \ast \Delta_1g,  \]
thus we need to show that $\Delta_1 g$ is a weakly positive
measure. Using Lemma \ref{propg} (a), we have
$g=h(L_1,\cdot)-h(L_2,\cdot)$ for some convex bodies $L_1,L_2 \in
\mathcal{K}^n$. Hence,
\[\Delta_1 g=S_1(L_1,\cdot)-S_1(L_2,\cdot).   \]
Using again Lemma \ref{delta1lem} and (\ref{delta1hk}), we obtain
\[S_1(\Phi K,\cdot)=S_{n-1}(K,\cdot) \ast \Delta_1 g \in \mathcal{M}^+_o(S^{n-1}).   \]
Thus, the desired result follows from Lemma \ref{weakpos3} and
from the fact that the set of surface area measures is a dense
subset of
$\mathcal{M}_o^+(S^{n-1})$. \hfill $\blacksquare$ \\

For $K, L \in \mathcal{K}^n$ and $i=0,\ldots,n-2$, let $W_i(K,L)$
denote the mixed volume $V(L,K,\ldots,K,B^n,\ldots,B^n)$, where
$K$ appears $n-1-i$ times and $B^n$ appears $i$ times. Note that,
$W_0(K,L)=V_1(K,L)$. For our further investigations we state the
following consequence of Lemma \ref{selfadlemma}.

\begin{lem} \label{durchsch} For
$i=0,\ldots, n-2$ and $K, L \in \mathcal{K}^n$,
\[W_i(K,\Phi_i L)=W_i(L,\Phi_i K).   \]
\end{lem}
{\it Proof}: Let $g \in
\mathcal{C}(S^{n-1},\mbox{\raisebox{-0.01cm}{$\down{e}$}})$
denote the generating function of $\Phi$. From
(\ref{mixedvolmixedsurf}), the definition of $W_i(K,L)$ and Lemma
\ref{selfadlemma}, it follows that
\begin{eqnarray*}
W_i(K,\Phi_i L)& = & \frac{1}{n}\langle h(\Phi_i L,\cdot),
S_{n-1-i}(K,\cdot) \rangle = \frac{1}{n} \langle
S_{n-1-i}(L,\cdot) \ast g,S_{n-1-i}(K,\cdot) \rangle \\
& = & \frac{1}{n} \langle S_{n-1-i}(L,\cdot) ,S_{n-1-i}(K,\cdot)
\ast g \rangle = W_i(L,\Phi_i K).
\end{eqnarray*}
\hfill $\blacksquare$ \\

The Shephard problem asks whether for $K, L \in \mathcal{K}_0^n$,
\begin{equation} \label{pikenthalteninpil}
\mbox{vol}_{n-1}(K |u^{\bot})=h(\Pi K,u) \leq h(\Pi L,u) =
\mbox{vol}_{n-1}(L|u ^{\bot})
\end{equation}
for every $u \in S^{n-1}$ implies $V(K) \leq V(L)$. Obviously,
(\ref{pikenthalteninpil}) is equivalent to $\Pi K \subseteq \Pi
L$. As was shown independently by Petty \cite{petty67} and
Schneider \cite{schneider67}, the answer to Shephard's problem is
no in general, but if the body $L$ is a zonoid, the answer is yes.
The crucial tool in the proof of the latter statement is a special
case of Lemma \ref{durchsch}. In fact, an analogous result can be
shown for general Blaschke-Minkowski homomorphisms.

\begin{koro} \label{pettyschneider} Let $K \in
\mathcal{K}^n_i$ and $L \in \Phi_i \mathcal{K}^n_i$. Then, for
$i=0,\ldots,n-2$,
\[\Phi_i K \subseteq \Phi_i L \quad \Rightarrow \quad W_i(K) \leq W_i(L)   \]
and $W_i(K)=W_i(L)$ only if $K$ and $L$ are translates.
\end{koro}
{\it Proof}: From the monotonicity of mixed volumes, Lemma
\ref{durchsch} and the fact that $L=\Phi_i L_0$ for some convex
body $L_0 \in \mathcal{K}_i^n$, it follows that
\[W_i(K,\Phi_i L_0) = W_i(L_0,\Phi_i K) \leq W_i(L_0,\Phi_i L)=W_i(L,\Phi_i L_0)=W_i(L).   \]
Using the generalized Minkowski inequality
\begin{equation} \label{genminkinequ}
W_i(K,L)^{n-i} \geq W_i(K)^{n-1-i}W_i(L),
\end{equation}
with equality if and only if $K$ and $L$ are homothetic, we thus
get
\[W_i(K) \leq W_i(L),   \]
with equality only if $K$ and $L$ are homothetic. But homothetic
bodies of equal $i$th quermassintegral must be translates of each
other. \hfill $\blacksquare$ \\

\noindent The special case $i=0, \Phi = \Pi$ of Corollary
\ref{pettyschneider} is the result of Schneider and Petty. The
following result is a generalization of Corollary
\ref{meanwidth}, which follows from (\ref{defphii}).

\begin{koro} \label{meanwidth2} For
$i=0,\ldots,n-2$ and $K \in \mathcal{K}^n$,
\[W_{n-1}(\Phi_i K)=r_{\Phi}W_{i+1}(K),  \]
where $r_{\Phi} \in \mathbb{R}^+$ is the radius of the ball $\Phi
B^n $.
\end{koro}

\noindent We will prove now an upper bound for the $i$th
quermassintegral of $\Phi_i K$.

\begin{satz} \label{upbound} For $i=0,\ldots, n-2$ and $K \in
\mathcal{K}^n$,
\[W_{i+1}(K)^{n-i} \geq \frac{\kappa_{n}^{n-1-i}}{r_{\Phi}^{n-i}}W_i(\Phi_i K),   \]
where $r_{\Phi} \in \mathbb{R}^+$ is the radius of the ball $\Phi
B^n $. There is equality only if $\Phi_iK$ is a ball.
\end{satz}
{\it Proof}: Let $K \in \mathcal{K}^n$ and $0 \leq i \leq n-2$.
From inequality (\ref{ineququermass}), we get by repeated
application, the inequality
\[W_{n-1}(K)^{n-i} \geq \kappa_n^{n-1-i}W_i(K),   \]
where, for $K \in \mathcal{K}^n_{n-1}$, there is equality if and
only if $K$ is a ball. Taking $K=\Phi_i K$ and using Corollary
\ref{meanwidth2}, gives the desired result. \hfill $\blacksquare$ \\

In the following we will investigate, for $K \in \mathcal{K}_i^n$,
the similarity invariant ratio
\[\psi_i(K)=\frac{W_i(\Phi_iK)}{W_i(K)^{n-1-i}}.   \]
By a standard technique, using Blaschke's selection theorem, it
can be shown that $\psi_i$ attains a minimum on $\mathcal{K}_i^n$.
From the next theorem follows that the extremal bodies $K$ of this
minimum have the property that $K$ and $\Phi_i^2 K$ are
homothetic.

\begin{satz} \label{psii} If $K \in \mathcal{K}_i^n$ and $0 \leq i \leq n-2$, then
\[\psi_i(K) \geq \psi_i(\Phi_i K),   \]
with equality if and only if $K$ and $\Phi_i^2 K$ are homothetic.
\end{satz}
{\it Proof}: Let $K, L \in \mathcal{K}_i^n$. From the generalized
Minkowski inequality (\ref{genminkinequ}) together with Lemma
\ref{durchsch}, we obtain
\[W_i(L,\Phi_i K)^{n-i}=W_i(K,\Phi_i L)^{n-i}\geq W_i(K)^{n-1-i}W_i(\Phi_i L),    \]
with equality if and only if $K$ and $\Phi_i L$ are homothetic.
Setting $L=\Phi_i K$, gives
\[W_i(\Phi_i K)^{n-i} \geq W_i(K)^{n-1-i}W_i(\Phi_i^2K),  \]
with equality if and only if $K$ and $\Phi_i^2 K$ are homothetic.
\hfill $\blacksquare$\\

In the case $i=n-2$, there is a result of Kiderlen
\cite{kiderlen05} that $K$ and $\Psi^2 K$, for a nontrivial weakly
monotone Minkowski endomorphism $\Psi$, are homothetic if and only
if $K$ is a ball, where the combinations of the identity map and
the reflection in the origin are the trivial Minkowski
endomorphisms. Thus, Theorem \ref{upbound} and Theorem \ref{psii}
together with Kiderlen's result imply Theorem \ref{minkinequ}.

In the proof of Theorem \ref{minkinequ}, we have used only that
$\Phi_{n-2}$ is a weakly monotone (nontrivial) Minkowski
endomorphism. In fact, inequality (\ref{sharpmin}) with equality
cases is valid for every nontrivial weakly monotone Minkowski
endomorphism, compare also \cite{schneider74}, p.70, for a related
result. The reason why we chose the more restrictive formulation
of Theorem \ref{minkinequ} is the author's belief that Petty's
conjectured projection inequality holds in a more general form
for every Blaschke-Minkowski homomorphism and its induced
operators
\begin{equation} \label{schconj}
W_{i+1}(K)^{n-i} \geq
\frac{\kappa_n^{n-1-i}}{r_{\Phi}^{n-i}}W_{i}(\Phi_{i}K) \geq
\kappa_nW_{i}(K)^{n-1-i},
\end{equation}
giving a family of strengthened versions of the classical
inequalities (\ref{ineququermass}) between consecutive
quermassintegrals. The inequalities (\ref{ineququermass}) are
special cases of (\ref{schconj}) for the Blaschke-Minkowski
homomorphisms $K \mapsto cB(W_1(K),o)$, where $c \in \mathbb{R}^+$
and $B(W_1(K),o)$ is the ball with center in the origin and radius
$W_1(K)$. It is possible to show that the case $i=0$ of
(\ref{schconj}) implies the inequality for all other values of
$i$, for a proof compare the argument for the projection body
operator in \cite{lutwak90b}, p.57.

\vspace{0.7cm}
\medskip\noindent{\bf Acknowledgements.} The work of
the author was supported by the Austrian Science Fund (FWF),
within the scope of the project "Affinely associated bodies",
Project Number: P16547-N12 and the project "Phenomena in high
dimensions" of the European Community, Contract Number:
MRTN-CT-2004-511953. For their helpful remarks the author is
obliged to Monika Ludwig and Rolf Schneider.

\hfill\parbox[t]{6truecm}{ Forschungsgruppe \hfill\par Konvexe
und Diskrete Geometrie \hfill\par Technische Universit\"at
Wien\hfill\par Wiedner Hauptstra\ss e 8--10/1046\hfill\par A--1040
Vienna, Austria\hfill\par franz.schuster@tuwien.ac.at\hfill}


\vspace{0.5cm}
\begin{thebibliography}{99}
\footnotesize{
\parskip-0.1cm{
\bibitem{berg69} C. Berg, \emph{Corps convexes et potentiels
sph\'eriques}, Danske Vid. Selsk. Mat.-Fys. Medd. {\bf 37}, 6
(1969), 1-64.

\bibitem{bourgainlindenstrauss88} J. Bourgain and J.
Lindenstrauss, \emph{Projection bodies}, Geometric aspects of
functional analysis (1986/87), Springer, Berlin (1988), 250-270.

\bibitem{bolker69} E.D. Bolker, \emph{A class of convex bodies},
Trans. Amer. Math. Soc. {\bf 145} (1969), 323-345.

\bibitem{bronshtein79} E.M. Bronshtein, \emph{Extremal H-convex
bodies (in Russian)}, Sibirskii Mat. Zh. {\bf 20} (1979),
412-415. English Translation: Siberian Math. J. {\bf 20}, 295-297.

\bibitem{gardner95} R.J. Gardner, \emph{Geometric Tomography},
Cambridge University Press, 1995.

\bibitem{goodeyschneider80} P. Goodey, R. Schneider, \emph{On the
intermediate area functions of convex bodies}, Math. Z. {\bf 173}
(1980), 185-194.

\bibitem{goodeyweil92} P. Goodey and W. Weil, \emph{The
determination of convex bodies from the mean of random sections},
Math. Proc. Camb. Phil. Soc. {\bf 112} (1992), 419-430.

\bibitem{goodeyweil93} P. Goodey and W. Weil, \emph{Zonoids and
generalizations}, Handbook of Convex Geometry (P.M. Gruber and
J.M. Wills, eds.), North-Holland, Amsterdam, 1993, 1297-1326.

\bibitem{grinbergzhang99} E. Grinberg and G. Zhang,
\emph{Convolutions, transforms, and convex bodies}, Proc. London
Math. Soc. (3) {\bf 78} (1999), 77-115.

\bibitem{groemer96} H. Groemer, \emph{Geometric Applications of Fourier
Series and Spherical Harmonics}, Cambridge University Press, 1996.

\bibitem{hadwiger51} H. Hadwiger, \emph{Vorlesungen \"uber Inhalt, Oberfl\"ache und Isoperimetrie}, Springer, Berlin, 1957.

\bibitem{schneiderhug02} D. Hug and R. Schneider, \emph{Stability
results involving surface area measures of convex bodies},
Rendiconti Del Circolo Matematica Di Palermo {\bf 70} (2002),
21-51.

\bibitem{kiderlen05} M. Kiderlen, \emph{Blaschke- and
Minkowski-Endomorphisms of convex bodies}, Trans. Amer. Math.
Soc., to appear.

\bibitem{klainrota97} D.A. Klain and G.C.Rota, \emph{Introduction
to geometric probability}, Cambridge University Press, Cambridge,
1997.

\bibitem{ludwig02} M. Ludwig, \emph{Projection bodies and
valuations}, Adv. Math. {\bf 172} (2002), 158-168.

\bibitem{ludwig04} M. Ludwig, \emph{Minkowski valuations}, Trans.
Amer. Math. Soc. {\bf 357} (2005), no. 10, 4191-4213.

\bibitem{lutwak85} E. Lutwak, \emph{Mixed projection inequalities}, Trans.
Amer. Math. Soc. {\bf 287} (1985), no. 1, 91-105.

\bibitem{lutwak90b} E. Lutwak, \emph{On quermassintegrals of mixed
projection bodies}, Geom. Dedicata {\bf 33} (1990), 51-58.

\bibitem{lutwak93} E. Lutwak, \emph{Inequalities for mixed projection bodies}, Trans.
Amer. Math. Soc. {\bf 339} (1993), no. 2, 901-916.

\bibitem{mcmullen80} P. McMullen, \emph{Continuous translation
invariant valuations on the space of compact convex sets}, Arch.
Math. {\bf 34}, 377-384.

\bibitem{mcmullen93} P. McMullen, \emph{Valuations and
dissections}, Handbook of Convex Geometry, Vol. B (P.M. Gruber and
J.M. Wills, eds.), North Holland, Amsterdam, 1993, 933-990.

\bibitem{mcmullenschneider83} P. McMullen and R. Schneider,
\emph{Valuations on convex bodies}, Convexity and its
applications (P.M. Gruber and J.M. Wills, eds.), Birkh\"auser,
1983, 170-247.

\bibitem{petty67} C.M. Petty, \emph{Projection bodies},
\emph{Proceedings, Coll. Convexity, Copenhagen, 1965}, Kobenhavns
Univ. Mat. Inst. (1967), 234-241.

\bibitem{petty72} C.M. Petty, \emph{Isoperimetric problems}, Proc.
Conf. on Convexity and Combinatorial Geometry, Univ. of Oklahoma,
June 1971 (1972), 26-41.

\bibitem{schneider67} R. Schneider, \emph{Zu einem Problem von
Shephard \"uber die Projektionen konvexer K\"orper}, Math. Z. {\bf
101} (1967), 71-82.

\bibitem{schneider74} R. Schneider, \emph{Equivariant
endomorphisms of the space of convex bodies}, Trans. Amer. Math.
Soc. {\bf 194} (1974), 53-78.

\bibitem{schneider74b} R. Schneider, \emph{Bewegungs\"aquivariante, additive und
stetige Transformationen konvexer Bereiche}, Arch. Math. {\bf 25}
(1974), 303-312.

\bibitem{schneider93} R. Schneider, \emph{Convex Bodies: The
Brunn-Minkowski Theory}, Cambridge University Press, 1993.

\bibitem{schuster05} F.E. Schuster, \emph{Volume Inequalities and Additive Maps of Convex Bodies}, in preparation.

\bibitem{weil74} W. Weil, \emph{Decomposition of convex bodies},
Mathematika {\bf 21} (1974), 19-25.

\bibitem{weil74b} W. Weil, \emph{\"Uber den Vektorraum der
Differenzen von St\"utzfunktionen konvexer K\"orper}, Math. Nachr.
{\bf 59} (1974), 353-369.

\bibitem{weil76} W. Weil, \emph{Kontinuierliche
Linearkombination von Strecken}, Math. Z. {\bf 148} (1976), 71-84.

\bibitem{weil80} W. Weil, \emph{On surface area measures of convex
bodies}, Geom. Dedicata {\bf 9} (1980), 299-306.

}}
\end{thebibliography}
\end{document}